\documentclass[review,onefignum,onetabnum]{siamart190516}

\usepackage{lipsum}
\usepackage{amsfonts}
\usepackage{graphicx}
\usepackage{epstopdf}
\usepackage{algorithmic}
\ifpdf
  \DeclareGraphicsExtensions{.eps,.pdf,.png,.jpg}
\else
  \DeclareGraphicsExtensions{.eps}
\fi

\newsiamremark{remark}{Remark}
\newsiamremark{hypothesis}{Hypothesis}
\crefname{hypothesis}{Hypothesis}{Hypotheses}
\newsiamthm{claim}{Claim}

\headers{On the Convergence of the Multipole Expansion Method}{B. Fitzpatrick, E. De Sena, and T. van Waterschoot}

\title{On the Convergence of the Multipole Expansion Method\thanks{Submitted to the editors DATE.
\funding{The research leading to these results has received funding from the European Research Council under the European Union's Horizon 2020 research and innovation program / ERC Consolidator Grant: SONORA (no. 773268).}}}

\author{Brian Fitzpatrick\thanks{Department of Electrical Engineering, KU Leuven, 3001 Leuven, Belgium (\email{brian.fitzpatrick@kuleuven.be}). Current address: School of Physics, TU Dublin, Grangegorman, Dublin 7, Ireland (\email{brian.fitzpatrick@tudublin.ie})}
\and Enzo De Sena \thanks{Institute of Sound Recording, University of Surrey, Guilford, GU2 7XH, UK (\email{e.desena@surrey.ac.uk}).}
\and Toon van Waterschoot \thanks{Department of Electrical Engineering, KU Leuven, 3001 Leuven, Belgium (\email{toon.vanwaterschootv@kuleuven.be}).}}

\usepackage{amsopn}

\usepackage[title]{appendix}
\usepackage{tabularx}
\usepackage[margin=0.5in]{geometry}
\usepackage{amssymb}
\usepackage{graphicx}
\usepackage{hyperref}
\usepackage{marvosym}

\usepackage{soul}
\usepackage{bbm}
\usepackage{bm}
\usepackage[english]{babel}
\usepackage{amsmath}
\usepackage{amssymb}
\usepackage{enumitem}
\usepackage{graphicx}
\usepackage{grffile}
\usepackage{color}
\usepackage{xcolor}
\usepackage{booktabs}
\usepackage{array}   
\usepackage{afterpage}
\usepackage{subfig}
\usepackage[font={small,it}]{caption}
\usepackage{dsfont}
\usepackage{pifont}
\usepackage{graphicx}
\usepackage{color}
\usepackage{xcolor}
\usepackage{mathtools}
\numberwithin{equation}{section}
\usepackage{hyperref}
\usepackage[framemethod=tikz]{mdframed}
\usepackage{lipsum}
\usepackage[most]{tcolorbox}

\newcommand{\ds}{\displaystyle}

\def \Vh0{\stackrel{\circ}{V}_h} \def\to{\rightarrow}

\def\cE{{\mathcal{E}}}

\def\cT{{\mathcal{T}}}

\newcommand{\bbA}{\mathbb{A}}
\newcommand{\bbB}{\mathbb{B}}
\newcommand{\bbC}{\mathbb{C}}

\newcommand{\bbI}{\mathbb{I}}
\newcommand{\bbN}{\mathbb{N}}
\newcommand{\bbR}{\mathbb{R}}
\newcommand{\bbW}{\mathbb{W}}

\newcommand{\bbZ}{\mathbb{Z}}

\newcommand{\N}{\mathbb{N}}

\newcommand{\beq}{\begin{equation}}
\newcommand{\eeq}{\end{equation}}
\newcommand{\bes}{\begin{equation*}}
\newcommand{\ees}{\end{equation*}}

\newcommand*{\blue}{\textcolor{blue}}

\newcommand*{\sumainb}[2]{\sum_{#1 \in #2}}

\newcommand*{\sumz}[1]{\sum_{#1\in \mathbb{Z}}}

\def\bf{{\mathbf{f}}}

\def\b0{{\mathbf{0}}}

\newcommand{\benumabc}{\begin{enumerate}[label=\textbf{(\alph*)}]}
\newcommand{\eenumabc}{\end{enumerate}}
\newcommand{\benumi}{\begin{enumerate}[label=\textbf{(\roman*)}]}
\newcommand{\eenumi}{\end{enumerate}}

\makeatletter
\def\BState{\State\hskip-\ALG@thistlm}
\makeatother

\newtcolorbox{lbluebox}[1][]{
    breakable,
    colback=blue!05,
    colframe=blue!05,
    outer arc=0pt,
    parbox=false,
    #1,
}

\newtcolorbox{bluebox}[1][]{
    breakable,
    colback=blue!20,
    colframe=blue!20,
    outer arc=0pt,
    parbox=false,
    #1,
}

\newtcolorbox{orangebox}[1][]{
    breakable,
    colback=orange!20,
    colframe=orange!20,
    outer arc=0pt,
    parbox=false,
    #1,
}

\newtcolorbox{redbox}[1][]{
    breakable,
    colback=red!20,
    colframe=red!20,
    outer arc=0pt,
    parbox=false,
    #1,
}

\newif\ifshowblue
\showbluefalse

\ifpdf
\hypersetup{
  pdftitle={On the Convergence of the Multipole Expansion Method},
  pdfauthor={B. Fitzpatrick, E. De Sena, and T. van Waterschoot}
}
\fi

\nolinenumbers

\begin{document}

\maketitle

\begin{abstract}
The multipole expansion method (MEM) is a spatial discretization technique that is widely used in applications that feature scattering of waves from circular cylinders. Moreover, it also serves as a key component in several other numerical methods in which scattering computations involving arbitrarily shaped objects are accelerated by enclosing the objects in artificial cylinders. A fundamental question is that of how fast the approximation error of the MEM converges to zero as the truncation number goes to infinity. Despite the fact that the MEM was introduced in 1913, and has been in widespread usage as a numerical technique since as far back as 1955, a precise characterization of the asymptotic rate of convergence of the MEM has not been obtained. In this work, we provide a resolution to this issue. While our focus in this paper is on the Dirichlet scattering problem, this is merely for convenience and our results actually establish convergence rates that hold for all MEM formulations irrespective of the specific boundary conditions or boundary integral equation solution representation chosen.
\end{abstract}

\begin{keywords}
Multiple scattering, multipole expansion, layer potential, error analysis, truncation error 
\end{keywords}

\begin{AMS}
31A10, 42B10, 65N12, 65N15, 65R20, 70F10, 78M15, 78M16
\end{AMS}

\section{Introduction} \label{sec:intro}
The multiple scattering of waves is an important topic that arises in a variety of scientific fields including acoustics, electromagnetics, elasticity, water waves, and quantum mechanics. In the frequency domain, it is well known that the scattering of waves from multiple disjoint circular cylinders (spheres in three dimensions) can be computed with exceptional efficiency using a meshless technique called the multipole expansion method (MEM), a spatial discretization technique that involves truncating infinite series of multipoles \cite{lai2019framework,barucq2018numerical,thierry2015mu,martin2006multiple,barucq2018numerical}.

The idea of applying multipole expansion techniques to multiple scattering problems can be traced back to 1913 with Z\.{a}vi\v{s}ka introducing it in \cite{zavivska1913beugung} to compute the scattering of waves from an array of parallel cylinders. In 1955, Row used the MEM to obtain a numerical solution to a multiple scattering problem \cite{row1955theoretical,martin2006multiple}. Since then, the MEM has appeared in countless fundamental and applied works which feature multiple scattering from cylinders or spheres.

The MEM also serves as a key component in several other numerical methods in which scattering computations involving arbitrarily shaped objects are accelerated by enclosing the objects in artificial cylinders or spheres. For instance, MEM-based formulations appear in scattering matrix methods \cite{lai2019framework,lai2015fast}, T-matrix methods \cite{ganesh2012convergence,cai1999large}, and Dirichlet-to-Neumann methods \cite{grote2004dirichlet,acosta2010coupling}.

Recently, the MEM has seen extensive use in the field of metamaterials where it has been combined with lattice summation techniques to allow for efficient computational simulations in problems featuring infinitely periodic lattices such as photonic and phononic crystals, and metasurfaces; see, for instance, \cite{ammari2017double} and the monograph \cite{ammari2018mathematical}. In the last few years, it has also been used to simulate scattering in the context of topological insulators \cite{ammari2019topologically,ammari2020robust}, and subwavelength resonance models of the cochlea \cite{ammari2020mimicking,ammari2019fully}. In addition, the MEM has been used in investigations of speckle statistics and non-invasive optical focusing in random scattering media \cite{bar2019monte,li2020non}. It is worth noting that these latter works utilized the open-source MEM scattering library $\mu$-diff which was released in 2015 \cite{thierry2015mu}. Another interesting application of $\mu$-diff can be found in \cite{haffner2021localization} where it was used to generate the training data for a neural network aimed at localizing multiple scattering objects. We use $\mu$-diff to validate our theoretical results in this paper; in fact our results establish the convergence theory for this library.

For an in-depth review of the MEM literature in the case of multiple cylinders, see \cite[Sec. 4.5.1]{martin2006multiple} and Sec. 4.5.1 of the associated 'Corrections and Additions' document for this monograph, since this was updated as recently as 2019 and features numerous additional examples of MEM usage in the literature.

The key parameter in any multipole-based scattering technique is the truncation number. The truncation number stipulates the number of terms that should be retained when one truncates the infinite series used to represent the problem in order to obtain a finite-dimensional discretization. Convergence theories for other multipole-based scattering methods such as the fast multipole method, the T-matrix method, and the method of fundamental solutions can  be found in \cite{darve2000fast,ganesh2012convergence,barnett2008stability}.

 In the case of the MEM, multiple scattering of waves between cylinders results in a coupling of the coefficients of the infinite series associated with each cylinder, and this phenomenon has a pronounced effect on the decay of the approximation error of the MEM as the truncation number $N \to \infty$. However, despite the fact that the MEM originated over a century ago, and has been in widespread use for well over fifty years, to the best of the authors' knowledge, the asymptotic rate of convergence of the approximation error of the MEM has not been properly quantified in the literature. Henceforth, we shall use the phrase \textit{'convergence of the MEM'} to refer to the rate of convergence of the approximation error of the MEM to zero.

Numerical investigations on the convergence of the MEM have recently been undertaken in \cite{barucq2018numerical,barucq2016study}. There have been several other works wherein numerical investigations have been performed to ascertain the performance of iterative methods applied to the MEM system of equations \cite{amirkulova2015acoustic,antoine2008numerical,lai2019framework}; it should be noted, however, that these iterative methods are converging to the approximate solution given by the MEM, hence, the underlying question of the asymptotic rate of convergence of the MEM solution to the true solution remains unanswered. In this paper we present a resolution to the long-standing problem of quantifying the decay of the MEM approximation error. The system of equations we consider first arose in Row's 1955 paper \cite[Eq. (3) and Eq. (5)]{row1955theoretical}. 

Let $\{\Omega_p\}_{p=1}^M$ be a set of disjoint circular cylinders in $\mathbb{R}^2$. Let the incident field be given by a point source located at $x_0$. Denote by $O_p$ the center of cylinder $\Omega_p$, and by $a_p$ its radius.  Denote by $d_{pq} = |O_p-O_q|$ the distance between the centers of cylinders $\Omega_p$ and $\Omega_q$. Denote by $d_{px_0} = |O_p-x_0|$ the distance between the center of cylinder $\Omega_p$ and the point source located at $x_0$. Denote by
\begin{align*}
\bbN_M & := \{n\in \bbN : 1 \le n \le M\}.
\end{align*}

Denote by $\cE(N)$ the approximation error of the MEM for the $M$ cylinder system. The following theorem provides an asymptotic bound on the convergence of $\cE(N)$.
\begin{theorem} \label{thm:mem-app-err-bnd}
As $N\to \infty$, it holds that
\begin{align} \label{eq:mem-app-err-bnd}
\cE(N) \lesssim \gamma_1(N):=
\begin{cases}
\ds \max \Bigg\{\max_{p\in \bbN_M} \bigg(\frac{a_p}{d_{px_0}}\bigg)^N, \  \max_{\substack{p,q\in\bbN_M \\ q \neq p}} \bigg(\frac{a_{p}}{d_{p q}-a_{q}}\bigg)^N\Bigg\}, \quad \quad & \text{for a point source}, \\
\ds \max_{\substack{p,q\in\bbN_M \\ q \neq p}} \bigg(\frac{a_{p}}{d_{pq}-a_{q}}\bigg)^N, \quad \quad & \text{for a plane wave}.
\end{cases}
\end{align}
\end{theorem}
The notation '$\lesssim$', which will be made precise later, can be interpreted as signifying the left hand side is less than the right hand side up to an asymptotically irrelevant sub-exponentially increasing factor. In the point source case, the first term on the right hand side of the inequality represents the approximation error directly associated with the incident field, while the second term represents the approximation error associated with the geometry; thus, the bound indicates that if the point source is sufficiently close to one of the cylinders, then it is this placement of the point source that ultimately dictates the rate of convergence, otherwise the convergence is dictated by the pair of cylinders that maximizes the expression $a_p/(d_{pq} - a_q)$. Numerical simulations verify that $\gamma_1(N)$ provides a tight characterization of the convergence of the MEM as long as the cylinders are not too close together.

When some cylinders are, in fact, in close proximity to one another, near-trapping of energy occurs as waves repeatedly reflect among theses cylinders, and thus it takes longer for energy to leak away to infinity. This behavior manifests itself in a decrease in the rate of convergence of the MEM. This decrease in the rate of convergence is captured by $\gamma_1(N)$, but having said that, it transpires that the bound becomes overly pessimistic. The issue is that, while the $a_p/(d_{pq} - a_q)$ term in $\gamma_1(N)$ characterizes interactions among the cylinders, it doesn't fully account for the phenomenon of multiple scattering. However, the derivation of \cref{thm:mem-app-err-bnd} relies on the asymptotic analysis of explicit expressions, and in the case of fully accounted for multiple scattering, analogous closed-form expressions are not available. To deal with this issue, and obtain a closed-form expression that yields a more accurate estimate of the convergence of the MEM in the case of closely spaced cylinders, we develop an approximation that accounts for first-order scattering effects while neglecting higher order multiple scattering effects. We then derive the rate of convergence of this approximation.

Denote by $\cE^{(1)}(N)$ the approximation error associated with the first-order scattering approximation of the MEM just discussed. Since first-order scattering effects strongly dominate over higher order scattering effects, we expect that $\cE^{(1)}(N)$ accurately characterizes the converge of the MEM, unless some of the cylinders are very close together. We have the following result for the convergence of $\cE^{(1)}(N)$.
\begin{theorem} \label{thm:mem-app-err-bnd-ps-pw}
As $N\to \infty$, it holds that
\begin{align} \label{eq:mem-app-err-bnd-ps-pw}
\cE^{(1)}(N) \lesssim \ds \gamma_2(N) := 
\begin{cases}
\ds \max \Bigg\{\max_{p\in \bbN_M} \bigg(\frac{a_p}{d_{px_0}}\bigg)^N, \max_{\substack{p,q\in\bbN_M \\ q \neq p}} \bigg(\frac{a_{p}d_{q x_0}}{d_{p q}d_{q x_0}-a_{q}^2}\bigg)^N\Bigg\}, \quad \quad & \text{for a point source}, \\
\ds \max_{\substack{p,q\in \bbN_M \\ q\neq p}} \bigg(\frac{a_p}{d_{pq}}\bigg)^{N}, \quad \quad & \text{for a plane wave}.
\end{cases}
\end{align}
\end{theorem}
Numerical simulations confirm that $\gamma_2(N)$ is indeed a far more accurate estimator of the convergence of the MEM than $\gamma_1(N)$ for the closely spaced case. When some cylinders are placed very close together, possibly almost touching, this approximation degrades somewhat as the higher multiple scattering effects become too significant to safely neglect. We discuss the possibility of accounting for higher multiple scattering effects by connecting our first-order scattering approximation with an iterative technique called the method of reflections \cite{ciaramella2017review,balabane2004boundary}.


Both \cref{thm:mem-app-err-bnd} and \cref{thm:mem-app-err-bnd-ps-pw} were derived for the case of a scattering problem featuring Dirichlet boundary conditions for which an indirect boundary integral equation solution representation was chosen. It is straightforward to show that if one were to consider a different set of boundary conditions, or a different boundary integral equation solution representation, different sub-exponential factors would arise during the derivations, but the expressions for $\gamma_1(N)$ and $\gamma_2(N)$ would ultimately remain the same. Thus, our theory holds not just for the setting we consider specifically, instead it holds for all boundary conditions and all boundary integral equation solution representations.

This paper is structured as follows. In \Cref{sec:prelim}, we introduce some notational conventions and function spaces. In \Cref{sec:prob-set-bie-form}, we describe the setting of the multiple scattering problem, and provide a representation of it in terms of an indirect boundary integral equation. In \Cref{sec:spatial-disc}, we apply a MEM discretization to the boundary integral equation. In \Cref{sec:mem-cvg-theory}, we present our convergence analysis of the MEM along with numerical simulations that validate the theoretical findings.

\section{Preliminaries} \label{sec:prelim}
Denote by
\begin{align*}
\bbZ_N := \{n\in \bbZ : -N \le n \le N \}, \quad \quad \bbZ_N^c := \bbZ\setminus \bbZ_N.
\end{align*}
Denote by $\Gamma_p:= \partial \Omega_p$ the boundary of cylinder $\Omega_p$, and by $\Gamma := \bigcup_{p=1}^M \Gamma_p$ the boundary of the full system of $M$ cylinders. Denote by $(r_p(x),\theta_p(x))$ the polar coordinates of the point $x \in \mathbb{R}^2$ with respect to a polar coordinate system whose origin is located at the center of cylinder $\Omega_p$; see \cite[Fig. 1]{thierry2015mu}.

The setting we use for the MEM formulation of a scattering problem is the space of periodic functions on $M$ cylinders, which is a natural generalization of the standard space of periodic functions on $[-\pi,\pi]$ \cite[Sec. 3.4 and Sec. 3.6]{iorio2001fourier}. Define the normalized basis functions $b_m^p$, for $p\in \bbN_M$ and $m\in \bbZ$, by
\begin{align} \label{eq:mem-basis-fns}
\begin{split}
b_m^p(x) :=
\begin{cases}
\dfrac{e^{im \theta_p(x)}}{\sqrt{2\pi a_p}}, \quad \quad & x \in \Gamma_p, \\
0, \quad \quad & x \in \Gamma \setminus \Gamma_p.
\end{cases}
\end{split}
\end{align}
Let $f$ be a periodic function that has an expansion in terms of these basis functions:
\begin{align*}
f(x) = \sum_{p\in \N_M} \sumz{m} f_m^p b_m^p(x), \quad \quad x \in \Gamma.
\end{align*}
The coefficient $f_m^p$ can be viewed as the $m$-th generalized Fourier coefficient associated with cylinder $\Omega_p$.
Now, letting $s \in \bbR$, we denote by $H^s(\Gamma)$ the fractional Sobolev space containing those $f$ that satisfy
\begin{align*}
||f||_s^2 & := \sum_{p\in \N_M} \sumz{m}(1+|m|^2)^s|f_m^p|^2 < \infty.
\end{align*}
\ifshowblue
\blue{[Double check normalization is ok]}
\fi
The inner product on $H^s(\Gamma)$ is 
\begin{align*}
(f, g)_s := \sum_{p\in \N_M} \sumz{m}(1+|m|^2)^s f_m^p \overline g_m^p,
\end{align*}
and the duality pairing on $H^s(\Gamma)\times H^{-s}(\Gamma)$ is
\begin{align*}
\langle f, g \rangle_{s,-s} = (f, g)_0 = \sum_{p\in \N_M} \sumz{m}  f_m^p \overline{g}_m^p, \quad \quad f \in H^s(\Gamma), \ g \in H^{-s}(\Gamma).
\end{align*}
We define $l_s$ to be the space consisting of doubly indexed sequences $(\alpha_m^p)_{m\in\bbZ}^{p\in \N_M}$ that satisfy
\begin{align} \label{eq:l-s-norm}
||(\alpha_m^p)_{m\in\bbZ}^{p\in \N_M}||_{l_s}^2 := \sum_{p\in \N_M} \sumz{m}(1+|m|^2)^s|\alpha_m^p|^2 < \infty.
\end{align}
Thus, we have that $f \in H^s(\Gamma)$ if and only if its associated Fourier coefficients $(f_m^p)_{m\in\bbZ}^{p\in \N_M} \in l_s$. As we primarly work in terms of Fourier coefficients in this paper, for convenience, we abuse notation and use the same notation for the norm of a function and the norm of its associated doubly indexed sequence of Fourier coefficients. Specifically, when we write $||(f_m^p)_{m\in\bbZ}^{p\in \N_M}||_s$, we mean $||(f_m^p)_{m\in\bbZ}^{p\in \N_M}||_{l_s}$. Likewise, for the operator norm of an operator $A:l_s \to l_t$, when we write $||A||_{s,t}$, we mean $||A||_{l_s \to l_t}$.

For our MEM approximation of functions $f \in H^s(\Gamma)$, we introduce the finite-dimensional spaces $\cT_N(\Gamma) \subset H^s(\Gamma)$ defined, for $N\ge 0$, as
\begin{align*}
\cT_N(\Gamma) := \bigg\{& f : f(x) = \sum_{p\in \N_M} \sum_{m\in \bbZ} f_m^p b_m^p(x), \ f_m^p \in \bbC \ \text{if} \ m\in \bbZ_N, \ f_m^p = 0 \ \text{if} \  m\in \bbZ_N^c, \ x \in \Gamma \bigg\}.
\end{align*}
We note the following identities where $J_m$ is the Bessel function of order $m$, and $H_m$ is the Hankel function of the first kind of order $m$:
\begin{align} \label{eq:spec-fns-neg-ord}
J_{-m}(x) = (-1)^m J_m(x), \quad \quad H_{-m}(x) = (-1)^m H_m(x).
\end{align}

\section{Problem setting and boundary integral equation formulation} \label{sec:prob-set-bie-form}
Recall that we are concerned with the scattering of waves by $M$ disjoint circular cylinders $\{\Omega_p\}_{p=1}^M$. Denote by $\Omega := \bigcup_{p=1}^M \Omega_p$ the set of all cylinders, and by $\Omega^+: = \bbR^2\setminus \overline \Omega$ the region exterior to the cylinders. Let $k >0$ be the wavenumber in $\Omega^+$. We assume that $k^2$ is not a Dirichlet eigenvalue of $-\Delta$ inside $\Omega_p$, for $p \in \bbN_M$; this condition ensures our problem is well-posed \cite[Prop. 2]{barucq2018numerical}. To derive the asymptotic rate of convergence of the MEM, it suffices to consider the Helmholtz equation with Dirichlet boundary conditions for the total field $u$:
\begin{alignat}{2} \label{eq:helmholtz}
\begin{split}
(\Delta + k^2) u = 0, & \quad \quad \text{in} \ \Omega^+, \\
u = 0, & \quad \quad \text{on} \ \Gamma, \\
u^\text{s} := u - u^\text{inc}, & \quad \quad \text{in} \ \overline{\Omega^+},
\end{split}
\end{alignat}
where $u^\text{s}$ is the scattered field, and $u^\text{inc}$ is the incident wavefield given by
\begin{align} \label{eq:u-inc-g}
u^{\text{inc}}(x) & :=
\begin{cases}
\ds e^{ik \beta \cdot x}, \quad \quad & \text{for a plane wave}, \\
\ds \frac{i}{4}H_0(k|x-x_0|), \quad \quad & \text{for a point source}.
\end{cases}
\end{align}
Here, $H_0$ is the Hankel function of the first kind of order zero, and $\beta = [\cos(\hat \beta),\sin(\hat \beta)]^T$, with $\hat \beta \in [0,2\pi]$, is the direction of propagation of the plane wave.
Finally, we require that the scattered field satisfies the Sommerfeld radiation condition:
\begin{align*}
\frac{\partial u^\text{s}}{\partial r} - iku^\text{s} = o(r^{-1/2}), \quad \quad \text{as} \ r := |x| \to \infty,
\end{align*}
where $\partial/\partial r$ is the radial derivative.
Denote by $G^k$ the outgoing fundamental solution of the associated Helmholtz equation:
\begin{equation} \label{eq:G}
G^k(x) := \frac{i}{4}H_0(k|x|).
\end{equation}
For $\phi \in H^{-1/2}(\Gamma)$, we introduce the single layer potential $S^k$ defined as
\begin{align} \label{eq:S-freq-dom}
(S^k\phi)(x) := & \int_{\Gamma} G^k(x-y) \ \phi(y) \ d\sigma(y), \quad \quad x \in \Omega^+.
\end{align}
Upon taking the trace of the single layer potential, we obtain the single layer boundary integral operator $V^k:H^{-1/2}(\Gamma) \to H^{1/2}(\Gamma)$ given by
\begin{align} \label{eq:V-freq-dom}
(V^k\phi)(x) := & \int_{\Gamma} G^k(x-y) \ \phi(y) \ d\sigma(y), \quad \quad x \in \Gamma.
\end{align}
We consider an indirect boundary integral equation representation for the solution of \cref{eq:helmholtz}:
\begin{align*}
u & = u^\text{inc} + u^\text{s} = u^\text{inc} + S^k \phi,\quad \quad \text{in} \ \Omega^+.
\end{align*}
Upon taking the trace of this equation, we obtain the following boundary integral equation:
\begin{align} \label{eq:bie}
V^k \phi & = f, \quad \quad \text{on} \ \Gamma,
\end{align}
where we denote by $f := - u^\text{inc}$ for convenience. Our aim is to ascertain a precise characterization of the asymptotic rate of convergence of the MEM applied to this equation. In the sequel we suppress the wavenumber dependence of $V^k$ for clarity and simply write $V$.

\section{Spatial discretization with the MEM} \label{sec:spatial-disc}
To obtain a MEM discretization, first we have to represent the scattering problem in terms of infinite series of multipoles. Once this representation has been obtained, we truncate the infinite series to obtain a finite-dimensional discretized problem that can be solved numerically. We obtain a weak formulation of the multiple scattering problem by multiplying \cref{eq:bie} by a test function and integrating over the boundary of the cylinders.

\vspace{\baselineskip}
\begin{quote}
Given $f \in H^{1/2}(\Gamma)$, find $\phi \in H^{-1/2}(\Gamma)$ such that
\begin{align} \label{eq:bie-weak-form}
\langle V \phi, \psi \rangle_{1/2,-1/2} = \langle f, \psi \rangle_{1/2,-1/2}, \quad \forall \ \psi \in H^{-1/2}(\Gamma).
\end{align}
\end{quote}
\vspace{\baselineskip}
\begin{align*}
\phi = \sum_{p\in \N_M} \sumz{m} \phi^p_m  b^p_m, \quad \psi = \sum_{p\in \N_M} \sumz{m} \psi^p_m  b^p_m, \quad f & = \sum_{p\in \N_M} \sumz{m} f^p_m  b^p_m,
\end{align*}
the problem becomes a matter of finding the doubly indexed sequence of coefficients $(\phi_m^p)_{m\in\bbZ}^{p\in \bbN_M}$. Specifically, substituting the above expansions into \cref{eq:bie-weak-form} and using Galerkin orthogonality the problem becomes:

\vspace{\baselineskip}
\begin{quote}
Given $(f_m^p)_{m\in\bbZ}^{p\in \bbN_M} \in l_{1/2}$, find $(\phi_n^q)_{n\in\bbZ}^{q\in \bbN_M} \in l_{-1/2}$ such that
\begin{align} \label{eq:bie-weak-form-coeffs}
\sum_{q\in \N_M} \sumz{n} \phi_n^q \langle V b_n^q, b_m^p \rangle_{1/2,-1/2} = f_m^p,
\end{align}
for $m \in\mathbb{Z}, \ p \in \bbN_M$.
\end{quote}
\vspace{\baselineskip}
Define the infinite dimensional per-cylinder coefficient vectors, for $p \in \bbN_M$, by
\begin{align*}
\phi^p & := [\dots,\phi_{-2}^p,\phi_{-1}^p,\phi_0^p,\phi_1^p,\phi_2^p,\dots]^T, \\
f^p & := [\dots,f_{-2}^p,f_{-1}^p,f_0^p,f_1^p,f_2^p,\dots]^T.
\end{align*} 
Then the problem can be expressed in terms of infinite block matrices and vectors:

\vspace{\baselineskip}
\begin{quote}
Given $F \in l_{1/2}$, find $\Phi \in l_{-1/2}$ such that
\begin{align} \label{eq:bie-weak-form-matrix}
\mathbb{V} \Phi = F.
\end{align}
\end{quote}
Here,
\begin{align*}
\mathbb{V} =
\begin{bmatrix}
V^{11} & V^{12} & \dots & V^{1M} \\
V^{21} & V^{22} & \dots & V^{2M} \\
\vdots & \vdots & \ddots & \vdots \\
V^{M1} & V^{M2} & \dots & V^{MM}
\end{bmatrix}, \quad
\Phi =
\begin{bmatrix}
\phi^1 \\
\phi^2 \\
\vdots \\
\phi^M \\
\end{bmatrix}, \quad
F =
\begin{bmatrix}
f^1 \\
f^2 \\
\vdots \\
f^M \\
\end{bmatrix}.
\end{align*}
The elements of the matrices $V^{pq}$ are given by
\begin{alignat*}{3}
V_{mn}^{pq} & := \langle V b_n^q, b_m^p \rangle_{1/2,-1/2}, \quad \quad m,n \in \bbZ.
\end{alignat*}
The elements $V_{mn}^{pq}$ have explicit representations \cite{thierry2015mu,barucq2018numerical}:
\begin{align} \label{eq:V-pq-mn}
V_{mn}^{pq} & =
\begin{cases}
\dfrac{i \pi a_p}{2} J_m(k a_p) H_{m}(k a_p), \quad \quad & p = q, \ m = n, \\
0, \quad \quad & p = q, \ m \neq n, \\
\dfrac{i \pi \sqrt{a_p a_q}}{2} \ds J_m(k a_p) H_{m-n}(kd_{pq}) e^{i(m-n)\theta_{pq}} J_n(k a_q), \quad \quad & p \neq q.
\end{cases}
\end{align}

Here, $\theta_{pq}$ is the angle between cylinder $\Omega_p$ and cylinder $\Omega_q$; see \cite[Figure 1]{thierry2015mu} or \cite[Figure 2.1]{martin2006multiple}. The incident field coefficients also have explicit representations \cite[Prop. 4 and Prop. 5]{thierry2015mu}:
\begin{align} \label{eq:f-p-m}
f^p_m & =
\begin{cases}
- \ds \sqrt{2\pi a_p} e^{-ik \beta \cdot O_p} e^{im (\pi/2 + \hat \beta)}J_m(ka_p), \quad \quad & \text{for a plane wave}, \\
- \ds \frac{i \pi a_p}{2}J_m(k a_p)H_m(k d_{px_0}) \ \overline{b_m^p(x_0)}, \quad \quad & \text{for a point source}.
\end{cases}
\end{align}
On a historical note, representations \cref{eq:V-pq-mn} and \cref{eq:f-p-m} for $V_{mn}^{pq}$ and $f_m^p$, respectively, appear as Equation (5) in Row's 1955 paper \cite{row1955theoretical}. Equation (3) in Row's paper corresponds to \cref{eq:bie-weak-form-matrix}.
It is common to apply a diagonal preconditioner to the system of equations \cref{eq:bie-weak-form-matrix}, as it vastly improves the performance of iterative solvers such as GMRES \cite{lai2019framework,thierry2015mu}. We apply the same preconditioner in this paper, however, our motivation for applying the preconditioner is rather to facilitate the convergence analysis in \Cref{sec:mem-cvg-theory}. In any case, we multiply both sides of \cref{eq:bie-weak-form-matrix} by
\begin{align} \label{eq:B-mat}
\bbB :=
\begin{bmatrix}
B^{11} & 0 & \dots & 0 \\
0 & B^{22} & \dots & 0 \\
\vdots & \vdots & \ddots & \vdots \\
0 & 0 & \dots & B^{MM}
\end{bmatrix}
:=
\begin{bmatrix}
(V^{11})^{-1} & 0 & \dots & 0 \\
0 & (V^{22})^{-1} & \dots & 0 \\
\vdots & \vdots & \ddots & \vdots \\
0 & 0 & \dots & (V^{MM})^{-1}
\end{bmatrix},
\end{align}
and obtain the following equivalent problem:

\vspace{\baselineskip}
\begin{quote}
Given $G \in \Phi \in l_{-1/2}$, find $\Phi \in l_{-1/2}$ such that
\begin{align} \label{eq:bie-weak-form-matrix-equiv}
\mathbb{W}\Phi = G.
\end{align}
\end{quote}
\vspace{\baselineskip}
Here, $G=\mathbb{B}F$, and $\mathbb{W} = \mathbb{B}\mathbb{V} = \mathbb{I} + \mathbb{A}$, with
\begin{align*}
\mathbb{I} =
\begin{bmatrix}
I^{11} & 0 & \dots & 0 \\
0 & I^{22} & \dots & 0 \\
\vdots & \vdots & \ddots & \vdots \\
0 & 0 & \dots & I^{MM}
\end{bmatrix}, \quad
\mathbb{A} =
\begin{bmatrix}
0 & A^{12} & \dots & A^{1M} \\
A^{21} & 0 & \dots & A^{2M} \\
\vdots & \vdots & \ddots & \vdots \\
A^{M1} & A^{M2} & \dots & 0
\end{bmatrix}, \quad
G =
\begin{bmatrix}
g^1 \\
g^2 \\
\vdots \\
g^M \\
\end{bmatrix}.
\end{align*}
The elements of the matrices $I^{pp}$ are given by
\begin{align*}
I_{mn}^{pp} & := 
\begin{cases}
1, & \quad \quad m = n, \\
0, & \quad \quad m \neq n.
\end{cases}
\end{align*}
The matrices $A^{pq}$ and vectors $g^p$ are given by
\begin{align}
A^{pq} & := B^{pp}V^{pq}, \label{eq:A-pq} \\
g^p & := B^{pp}f^p. \label{eq:g-p}
\end{align}
Therefore, by \cref{eq:g-p}, \cref{eq:B-mat}, and \cref{eq:f-p-m} we have:
\begin{align} \label{eq:g-m-p}
g_m^p = &
\begin{cases}
- \ds \frac{2\sqrt{2}}{i \sqrt{\pi a_p}} e^{-ik \beta \cdot O_p} e^{im (\pi/2 + \hat \beta)}(H_m(ka_p))^{-1} , \quad \quad & \text{for a plane wave}, \\
- \ds \frac{H_m(k d_{px_0})}{H_m(ka_p)} \ \overline{b_m^p(x_0)}, \quad \quad & \text{for a point source}.
\end{cases}
\end{align}
Likewise, by \cref{eq:A-pq}, \cref{eq:B-mat}, and \cref{eq:V-pq-mn} we have:
\begin{align} \label{eq:A-pq-mn}
A_{mn}^{pq} & =
\begin{cases}
0, \quad \quad & p = q, \\
\sqrt{\dfrac{a_p}{a_q}} \ds (H_m(k a_p))^{-1} H_{m-n}(kd_{pq}) e^{i(m-n)\theta_{pq}} J_n(k a_q), \quad \quad & p \neq q.
\end{cases}
\end{align}
So far we have just expressed the continuous problem \cref{eq:bie-weak-form} in a different form. The next stage in the MEM discretization procedure consists of truncating the infinite-dimensional block matrices and block vectors to obtain a finite-dimensional discretized problem. However, in order to perform a convergence analysis, rather than directly working with the finite-dimensional truncated objects, we use infinite-dimensional versions of them in which the elements that fall outside the truncation range are set to $0$. Throughout this paper, we use tildes to denote the effectively finite-dimensional truncated MEM matrices and vectors associated with the infinite-dimensional matrices and vectors of the original problem.

We denote by $N \ge 0$ the MEM truncation number. Define
\begin{align*}
\tilde \phi^p & := [\dots,0,0,\phi_{-N}^p,\dots,\phi_{-1}^p,\phi_0^p,\phi_1^p,\dots,\phi_N^p,0,0,\dots]^T, \\
\tilde g^p & := [\dots,0,0,g_{-N}^p,\dots,g_{-1}^p,g_0^p,g_1^p,\dots,g_N^p,0,0,\dots]^T.
\end{align*}
Note that the various mathematical objects in the discrete problem have a dependence on $N$ but we regularly suppress this in the sequel. Truncation of the matrices and vectors that comprise the block matrices and block vectors in \cref{eq:bie-weak-form-matrix-equiv} leads to the following discrete problem:

\vspace{\baselineskip}
\begin{quote}
Given $\tilde G \in \cT_N(\Gamma)$, find $\tilde \Phi \in \cT_N(\Gamma)$ such that
\begin{align} \label{eq:bie-weak-form-matrix-equiv-disc}
\tilde{\mathbb{W}} \tilde \Phi = \tilde G.
\end{align}
\end{quote}
\vspace{\baselineskip}
Here, $\tilde{\mathbb{W}} = \tilde{\mathbb{I}} + \tilde{\mathbb{A}}$, with
\begin{align*}
\tilde{\mathbb{I}} =
\begin{bmatrix}
\tilde I^{11} & 0 & \dots & 0 \\
0 & \tilde I^{22} & \dots & 0 \\
\vdots & \vdots & \ddots & \vdots \\
0 & 0 & \dots & \tilde I^{MM}
\end{bmatrix}, \quad
\tilde{\bbA} =
\begin{bmatrix}
0 & \tilde A^{12} & \dots & \tilde A^{1M} \\
\tilde A^{21} & 0 & \dots & \tilde A^{2M} \\
\vdots & \vdots & \ddots & \vdots \\
\tilde A^{M1} & \tilde A^{M2} & \dots & 0
\end{bmatrix}, \quad
\tilde G & =
\begin{bmatrix}
\tilde g^1 \\
\tilde g^2 \\
\vdots \\
\tilde g^M \\
\end{bmatrix}.
\end{align*}
The elements of the matrices $\tilde I^{pp}$ are given by
\begin{align*}
\tilde I_{mn}^{pp} & := 
\begin{cases}
1, & \quad \quad m, n \in \mathbb{Z}_N, \ m = n, \\
0, & \quad \quad m, n \in \mathbb{Z}_N, \ m \neq n, \\
0, & \quad \quad m \in \mathbb{Z}_N^c, \ \text{or} \ n \in \mathbb{Z}_N^c.
\end{cases}
\end{align*}
The elements of the matrices $\tilde A^{pq}$ are given by
\begin{align} \label{eq:A-pq-mn-N}
\tilde A_{mn}^{pq} & := 
\begin{cases}
A_{mn}^{pq}, & \quad \quad m,n \in \mathbb{Z}_N, \\
0, & \quad \quad m \in \mathbb{Z}_N^c, \ \text{or} \ n \in \mathbb{Z}_N^c.
\end{cases}
\end{align}
The elements of the vectors $\tilde g^p$ are given by
\begin{align} \label{eq:g-p-m-N}
\tilde g_m^p & := 
\begin{cases}
g_m^p, & \quad \quad m \in \mathbb{Z}_N, \\
0, & \quad \quad m \in \mathbb{Z}_N^c.
\end{cases}
\end{align}

In practise, one numerically solves the linear system of equations in \cref{eq:bie-weak-form-matrix-equiv-disc} to obtain the MEM approximate solution $\tilde \Phi(N)$. The approximation error $\cE(N)$ is the difference between the solution of the original problem $\Phi$ and the MEM approximate solution $\tilde \Phi(N)$:
\begin{align} \label{eq:app-err}
\cE(N) := ||\Phi - \tilde \Phi(N)||_{-1/2}.
\end{align}
The question of precisely how fast $\cE(N)$ decays to zero as $N\to \infty$ is a fundamental aspect of the MEM that has not been properly addressed to date in the literature and is the focus of the next section.

\section{MEM convergence theory} \label{sec:mem-cvg-theory}
Before we begin, we introduce some notation for the purposes of clarity. The functions involved in our MEM convergence theory frequently feature rates of growth or decay that are at least exponential with respect to some variable of interest. In light of this, algebraic factors in these functions are asymptotically irrelevant; they ultimately lead to arbitrarily small corrections to the rate of convergence. This motivates the introduction of the following notational convention, since it allows us to absorb algebraic factors. We use the notation
\begin{align} \label{eq:exp-bnd-notation}
a \lesssim b,
\end{align}
to signify that $a$ is bounded by $b$ up to some function that increases sub-exponentially with respect to $m$. To be specific, we define a sub-exponentially increasing function as any function that increases with respect to $m$ slower than $e^{Cm}$, for $C>0$, as $m \to \infty$.

As an example of this notation let us consider the large order asymptotics of the Bessel function $J_m$ and Hankel function $H_m$ of order $m$ as we will require these later. For $m \in \bbZ\setminus \{0\}$ and $r>0$, the following super-exponential uniform bounds hold
\begin{alignat*}{3}
& c_J \frac{1}{\sqrt{|m|}}\bigg(\frac{e r}{2|m|}\bigg)^{|m|} \le && |J_m(r)| && \le C_J \frac{1}{\sqrt{|m|}}\bigg(\frac{e r}{2|m|}\bigg)^{|m|}, \\
& c_H \frac{1}{\sqrt{|m|}}\bigg(\frac{2|m|}{e r}\bigg)^{|m|} \le && |H_m(r)| && \le C_H \frac{1}{\sqrt{|m|}}\bigg(\frac{2|m|}{e r}\bigg)^{|m|},
\end{alignat*}
for some constants $c_J,c_H,C_J,C_H$; see \cite[Section 9.3]{abramowitz1948handbook} or \cite{barnett2008stability}. The notation of \cref{eq:exp-bnd-notation} allows us to 'disregard' the algebraic factors and say
\begin{equation} \label{eq:exp-bnds-J-H}
\begin{alignedat}{3}
\bigg(\frac{e r}{2|m|}\bigg)^{|m|} \lesssim \ |J_m(r)| \lesssim \bigg(\frac{e r}{2|m|}\bigg)^{|m|}, \quad \quad
\bigg(\frac{2|m|}{e r}\bigg)^{|m|} \lesssim \ |H_m(r)| \lesssim \bigg(\frac{2|m|}{e r}\bigg)^{|m|}.
\end{alignedat}
\end{equation}
Stirling's approximation will also be required later \cite{robbins1955remark}. Stirling's approximation states that for $m>0$, it holds that
\begin{align*}
\sqrt{2\pi} m^{m+1/2} e^{-m} \le m! \le e m^{m+1/2}e^{-m}.
\end{align*}
Using the notation introduced above, we can write this as
\begin{align} \label{eq:stir-approx}
m^m e^{-m} \lesssim m! \lesssim m^m e^{-m}.
\end{align}

\subsection{Bounding the approximation error of the MEM}
Recall the equations for the original problem \cref{eq:bie-weak-form-matrix-equiv}, and the discretized problem \cref{eq:bie-weak-form-matrix-equiv-disc}:
\begin{align*}
\bbW \Phi = (\bbI + \bbA)\Phi = G, \quad \quad \tilde{\bbW} \tilde \Phi = (\tilde{\bbI} + \tilde{\bbA})\tilde \Phi = \tilde G.
\end{align*}
Using these relations it is straightforward to show that $(\bbI + \tilde{\bbA})(\Phi - \tilde{\Phi}) = G -\tilde G + (\tilde{\bbA} - \bbA)\Phi$,
and thus
\begin{align*}
(\Phi - \tilde{\Phi}) = (\bbI + \tilde{\bbA})^{-1}(G -\tilde G + (\tilde{\bbA} - \bbA)\Phi).
\end{align*}
Upon taking norms and applying the triangle inequality, we obtain
\begin{align} \label{eq:app-err-bnd-0}
\cE(N) \le ||(\bbI + \tilde{\bbA}(N))^{-1}||_{-1/2,-1/2} \ \bigg(||G -\tilde G(N)||_{-1/2} + ||(\bbA-\tilde{\bbA}(N))\Phi||_{-1/2} \bigg).
\end{align}
We can also immediately say
\begin{align} \label{eq:app-err-bnd-1}
\cE(N) \le ||(\bbI + \tilde{\bbA}(N))^{-1}||_{-1/2,-1/2} \ \bigg(||G -\tilde G(N)||_{-1/2} + ||\bbA-\tilde{\bbA}(N)||_{-1/2,-1/2} \ ||\Phi||_{-1/2}\bigg).
\end{align}
Our plan in what follows is to derive the asymptotic bound for $\cE(N)$ given in  \cref{thm:mem-app-err-bnd} by explicitly estimating the right hand side of \cref{eq:app-err-bnd-1} as $N\to \infty$. It transpires that this bound provides a tight characterization of the approximation error when the cylinders are not too close together, however, it becomes somewhat pessimistic when some cylinders are, in fact, in close proximity to one another. Thus, once we have obtained an asymptotic bound on $\cE(N)$ using \cref{eq:app-err-bnd-1}, we will return to the closely spaced case and consider a first-order scattering approximation based on \cref{eq:app-err-bnd-0}; this approximation allows us to derive the bound given in \cref{thm:mem-app-err-bnd-ps-pw} which provides a more accurate representation of the convergence of the MEM in this particular regime.
\begin{proof}[Proof of \cref{thm:mem-app-err-bnd}]
It is well-known that when $k^2$ is not an interior Dirichlet eigenvalue for $-\Delta$ inside $\Omega_p$, for $p \in \mathbb{N}_M$, the operator $\bbA$
is compact due to the positive distance between any two cylinders in the system, and hence the fact that $\bbI + \bbA$ is invertible with bounded inverse follows by the Fredholm theory; see \cite[Prop. 2]{barucq2018numerical}, \cite[Sec. 5.2]{lai2019framework} or \cite[Thm. 2]{acosta2015surface}. Recall from \cref{eq:A-pq-mn-N} that $\tilde \bbA(N)$ is simply $\bbA$ with the elements of its sub-matrices set to zero outside a finite range. Hence, for sufficiently large $N$, $\bbI + \tilde \bbA(N)$ is also invertible and we have
\begin{align} \label{eq:I-m-A-inv-bnd}
||(\bbI + \tilde{\bbA}(N))^{-1}||_{-1/2,-1/2} & \le C, \quad \quad N \to \infty,
\end{align}
for some positive constant $C$. In \cref{lem:G-minus-G-N} and \cref{lem:A-minus-A-N}, we will establish that
\begin{align}
& ||G -\tilde G(N)||_{-1/2} \lesssim
\begin{cases}
\ds \max_{p\in \bbN_M} \bigg(\frac{eka_p}{2N}\bigg)^N, \quad \quad & \text{for a plane wave}, \\
\ds \max_{p\in \bbN_M} \bigg(\frac{a_p}{d_{px_0}}\bigg)^N, \quad \quad & \text{for a point source}, \label{eq:G-m-G-N-bnd}
\end{cases} \\
& ||\bbA-\tilde{\bbA}(N)||_{-1/2,-1/2} \lesssim \bigg(\frac{a_{p}}{d_{pq}-a_{q}}\bigg)^N. \label{eq:A-m-A-N-bnd}
\end{align}
In the case of a plane wave incident field, the bound on $||G -\tilde G(N)||_{-1/2}$ decays super-exponentially as $N\to \infty$, and therefore this bound will always be dominated by the bound on $||\bbA-\tilde{\bbA}(N)||_{-1/2,-1/2}$ which decays merely exponentially.
For an incident field due to a point source on the other hand, if the point source is sufficiently close to one of the cylinders, the bound on $||G -\tilde G(N)||_{-1/2}$ can be larger than $||\bbA-\tilde{\bbA}(N)||_{-1/2}$ as $N\to \infty$, and thus it can't be neglected. Bearing this in mind, and substituting \cref{eq:I-m-A-inv-bnd}, \cref{eq:G-m-G-N-bnd}, and \cref{eq:A-m-A-N-bnd} into \cref{eq:app-err-bnd-1} gives \cref{eq:mem-app-err-bnd}. 
\end{proof}

\begin{lemma} \label{lem:G-minus-G-N}
As $N\to \infty$, it holds that
\begin{align} \label{eq:lem:G-minus-G-N}
||G -\tilde G(N)||_{-1/2} \lesssim
\begin{cases}
\ds \max_{p\in \bbN_M} \bigg(\frac{eka_p}{2N}\bigg)^N, \quad \quad & \text{for a plane wave}, \\
\ds \max_{p\in \bbN_M} \bigg(\frac{a_p}{d_{px_0}}\bigg)^N, \quad \quad & \text{for a point source}.
\end{cases}
\end{align}
\begin{proof}
First, note that by the definition of $g_m^p$ in \cref{eq:g-m-p} and the uniform bounds in \cref{eq:exp-bnds-J-H}, for $m \ge 1$, we have that for a plane wave,
\begin{align} \label{eq:g-m-p-pw-bnd}
|g_m^p| & \lesssim \bigg(\frac{eka_p}{2m}\bigg)^m,
\end{align}
and for a point source,
\begin{align} \label{eq:g-m-p-ps-bnd}
|g_m^p| & \lesssim \bigg(\frac{2m}{ekd_{px_0}}\bigg)^m \bigg(\frac{eka_p}{2m}\bigg)^m = \bigg(\frac{a_p}{d_{px_0}}\bigg)^m.
\end{align}
Recalling \cref{eq:g-p-m-N}, we have
\begin{align*}
||G -\tilde G(N)||_{-1/2}^2 = \sum_{p\in\bbN_M}\sum_{m\in\bbZ} (1+|m|^2)^{-1/2} |g_m^p - \tilde g_m^p|^2
= \sum_{p\in\bbN_M}\sum_{m > N} (1+|m|^2)^{-1/2} (|g_{-m}^p|^2 + |g_{m}^p|^2).
\end{align*}
Now, by the definition of $|g_m^p|$ in \cref{eq:g-m-p}, and using the relations in \cref{eq:spec-fns-neg-ord},
\begin{align*} \label{eq:G-m-G-N-interim-bnd}
||G -\tilde G(N)||_{-1/2}^2 & = \sum_{p\in\bbN_M}\sum_{m > N} (1+|m|^2)^{-1/2} |g_m^p|^2 \le \max_{p\in \bbN_M} \sum_{m > N} (1+|m|^2)^{-1/2} |g_m^p|^2.
\end{align*}
In this case of a plane wave, by \cref{eq:g-m-p-pw-bnd}, this means
\begin{align*}
||G -\tilde G(N)||_{-1/2}^2 & \lesssim \max_{p\in \bbN_M} \sum_{m > N} (1+|m|^2)^{-1/2} \bigg(\frac{eka_p}{2m}\bigg)^{2m} \le \max_{p\in \bbN_M} C (1+N^2)^{-1/2} \bigg(\frac{eka_p}{2N}\bigg)^{2N},
\end{align*}
for some sufficiently large constant $C$, as $N \to \infty$. Then, upon absorbing the algebraic factor using \cref{eq:exp-bnd-notation}, we obtain the result in \cref{eq:lem:G-minus-G-N} for the case of a plane wave. The result for the point source case can be obtained in a 
 fashion, albeit using \cref{eq:g-m-p-ps-bnd} instead of \cref{eq:g-m-p-pw-bnd}.
\end{proof}
\end{lemma}
Before proving \cref{lem:A-minus-A-N}, we need the following two lemmas.
\begin{lemma} \label{lem:A-pq-mn-series-bnd}
For $p \neq q$, it holds that
\begin{align} \label{eq:h-1-pq-N-pre_bound}
\sumainb{m}{\bbZ_N^c} \sumz{n} |A_{mn}^{pq}|^2 & \lesssim \sum_{m>N} \bigg(\frac{a_p}{d_{pq}}\bigg)^{2m} + \sum_{m>N} \sum_{n=1}^\infty \sigma^{pq}(m,n),
\end{align}
where
\begin{align} \label{eq:sigma-pq}
\sigma^{pq}(m,n) = \bigg(\frac{m+n}{m}\bigg)^{2m}\bigg(\frac{m+n}{n}\bigg)^{2n}\bigg(\frac{a_p}{d_{pq}}\bigg)^{2m}\bigg(\frac{a_q}{d_{pq}}\bigg)^{2n},
\end{align}
and $A_{mn}^{pq}$ is given in \cref{eq:A-pq-mn}.

\begin{proof}
First, define $C_{pq} := \pi^2 a_p a_q/4$. Then, using the relations in \cref{eq:spec-fns-neg-ord}, and the fact that
\begin{align} \label{eq:abs_h_m_str_inc}
H_{m-n}(x) \le H_{m+n}(x), \quad \quad m, n \ge 0,
\end{align}
it is straightforward to show that
\begin{align*}
\sumainb{m}{\bbZ_N^c} \sumz{n} |A_{mn}^{pq}|^2 = & \ C_{pq} \sumainb{m}{\bbZ_N^c} \sumz{n} \bigg|\frac{H_{m-n}(k d_{pq})J_n(ka_q)}{H_m(ka_p)}\bigg|^2 \\
\le & \ 2 C_{pq} \bigg(\sum_{m>N} \bigg|\frac{H_{m}(k d_{pq}) J_0(ka_q)}{H_m(ka_p)}\bigg|^2 + \sum_{m>N} \sum_{n=1}^\infty \bigg|\frac{H_{m+n}(k d_{pq}) J_n(ka_q)}{H_m(ka_p)}\bigg|^2\bigg).
\end{align*}
Then, upon applying the uniform bounds for the Bessel and Hankel functions \cref{eq:exp-bnds-J-H}, and absorbing constant factors using \cref{eq:exp-bnd-notation}, we have that for $m,n >0$,
\begin{align*}
2 C_{pq} \bigg|\frac{H_{m}(k d_{pq})}{H_m(ka_p)}J_0(ka_q)\bigg|^2 & \lesssim \bigg(\frac{2m}{ekd_{pq}}\bigg)^{2m} \bigg(\frac{ek a_p}{2m}\bigg)^{2m} = \bigg(\frac{a_p}{d_{pq}}\bigg)^{2m},
\end{align*}
and
\begin{align*}
2 C_{pq} \bigg|\frac{H_{m+n}(k d_{pq})}{H_m(ka_p)}J_n(ka_q)\bigg|^2
& \lesssim \bigg(\frac{2(m+n)}{ekd_{pq}}\bigg)^{2(m+n)} \bigg(\frac{ek a_p}{2m}\bigg)^{2m} \bigg(\frac{ek a_q}{2n}\bigg)^{2n} \\
& = \bigg(\frac{m+n}{m}\bigg)^{2m} \bigg(\frac{m+n}{n}\bigg)^{2n}\bigg(\frac{a_p}{d_{pq}}\bigg)^{2m}\bigg(\frac{a_q}{d_{pq}}\bigg)^{2n}.
\end{align*}
\end{proof}
\end{lemma}

In the next lemma, we derive an explicit representation of $\sum_{n=1}^\infty \sigma^{pq}(m,n)$ from \cref{lem:A-pq-mn-series-bnd}, as $m\to \infty$.
\begin{lemma} \label{lem:sigma-series-bnd}
As $m\to \infty$ it holds that
\begin{align*}
\sum_{n=1}^\infty \sigma^{pq}(m,n) \lesssim \bigg(\frac{a_p}{d_{pq}-a_q}\bigg)^{2m}.
\end{align*}
\begin{proof}
First, we rewrite $\sigma^{pq}$ as
\begin{align*}
\sigma^{pq}(m,n) & = (m+n)^{2(m+n)} \bigg(\frac{1}{m}\bigg)^{2m}\bigg(\frac{1}{n}\bigg)^{2n}\bigg(\frac{a_p}{d_{pq}}\bigg)^{2m}\bigg(\frac{a_q}{d_{pq}}\bigg)^{2n}.
\end{align*}
Stirling's approximation \cref{eq:stir-approx} gives
\begin{align*}
(m+n)^{2(m+n)} \lesssim ((m+n)!)^2 e^{2(m+n)}, \quad \quad \bigg(\frac{1}{m}\bigg)^{2m} \lesssim \frac{1}{(m!)^2 e^{2m}}, \quad \quad \bigg(\frac{1}{n}\bigg)^{2n} \lesssim \frac{1}{(n!)^2 e^{2n}}.
\end{align*}
Then, denoting by $z^{pq}:= (a_q/d_{pq})^2$, we get
\begin{align*}
\sum_{n=1}^\infty \sigma^{pq}(m,n)
& \lesssim \sum_{n=1}^\infty ((m+n)!)^2 e^{2(m+n)} \frac{1}{(m!)^2 e^{2m}} \frac{1}{(n!)^2 e^{2n}} \bigg(\frac{a_p}{d_{pq}}\bigg)^{2m} (z^{pq})^{n} \\
& = \bigg(\frac{a_p}{d_{pq}}\bigg)^{2m} \sum_{n=1}^\infty \bigg(\frac{(m+n)!}{m!n!}\bigg)^2 (z^{pq})^{n}.
\end{align*}
To find a simplified expression for the series in the above expression, consider that
\begin{align*}
\sum_{n=1}^\infty \bigg(\frac{(m+n)!}{m!n!}\bigg)^2 (z^{pq})^n = \frac{1}{(m!)^2} \sum_{n=1}^\infty \frac{((m+n)!)^2}{n!} \frac{(z^{pq})^n}{n!}
 = \sum_{n=1}^\infty \frac{(m+1)_n(m+1)_n}{(1)_n} \frac{(z^{pq})^n}{n!},
\end{align*}
where $(a)_n$ is the Pochhammer symbol, otherwise known as the rising factorial.
This is essentially the definition of the hypergeometric function ${}_2F_1$ for a particular set of parameters:
\begin{align*}
{}_2F_1(m+1,m+1;1;z^{pq}) = 1 + \sum_{n=1}^\infty \frac{(m+1)_n(m+1)_n}{(1)_n} \frac{(z^{pq})^n}{n!}.
\end{align*}
Hence we have
\begin{align*}
\sum_{n=1}^\infty \sigma^{pq}(m,n) \lesssim \bigg(\frac{a_p}{d_{pq}}\bigg)^{2m} {}_2F_1(m+1,m+1;1;(a_q/d_{pq})^2).
\end{align*}
To progress further we need a large argument asymptotic approximation of the hypergeometric function ${}_2F_1$ for the above set of parameters. This can be found in \cref{lem:hyper-geom-asy-ps} in the appendix, where we prove that as $m\to \infty$, it holds that
\begin{align*}
{}_{2}F_1(m+1,m+1;1;(a/b)^2) \lesssim \bigg(\frac{b^2}{b^2-a^2}\bigg)^m \bigg(\frac{b+a}{b-a}\bigg)^m.
\end{align*}
Therefore, as $m\to \infty$, we find that
\begin{align*}
\sum_{n=1}^\infty \sigma^{pq}(m,n) & \lesssim \bigg(\frac{a_p}{d_{pq}}\bigg)^{2m} \bigg(\frac{d_{pq}^2}{d_{pq}^2-a_q^2}\bigg)^m \bigg(\frac{d_{pq}+a_q}{d_{pq}-a_q}\bigg)^m = \bigg(\frac{a_p}{d_{pq}-a_q}\bigg)^{2m}.
\end{align*}
\end{proof}
\end{lemma}

\begin{remark}
It is worth nothing that if one were to choose a different set of boundary conditions than the Dirichlet conditions specified in \cref{eq:helmholtz}, or choose a different boundary integral equation solution representation for the problem, the result in \cref{lem:A-pq-mn-series-bnd} would not change. This is due to the fact that when we use the large order asymptotics of the Bessel and Hankel functions, the expressions obtained with those other choices differ from the expressions in \cref{lem:A-pq-mn-series-bnd} only by asymptotically irrelevant sub-exponentially increasing factors which are ultimately absorbed using the notation \cref{eq:exp-bnd-notation}. It is for this reason that our convergence theory holds for all MEM formulations irrespective of the specific boundary conditions or boundary integral equation solution representation chosen.
On a related note, recall that the specific boundary integral formulation we work with in this paper is an indirect single layer potential solution representation with a diagonal preconditioner applied. Interestingly, it has been shown \cite{thierry2014remark} that every boundary integral equation shares the same spectral properties after applying the relevant diagonal preconditioner, which results in iterative Krylov subspace solvers such as GMRES having the same rate of convergence, irrespective of the particular boundary integral equation formulation chosen.
\end{remark}
\begin{remark}
Results somewhat similar to the ones we derived in \cref{lem:A-pq-mn-series-bnd} and \cref{lem:sigma-series-bnd} were obtained in \cite{antoine2013dense} in the context of spectral and condition number estimates of the single layer operator in dense media in the low-frequency regime, using an approach that doesn't involve hypergeometric functions. The fact that the low-frequency results in \cite{antoine2013dense} are similar to our results based on large order asymptotics is not so surprising when one considers that $(i\pi a_p/2) J_m(ka_p) H_m(ka_p) \sim a_p/2m$, when either $k \to 0$ while $m > 0$, or $k > 0$ while $m \to \infty$, that is, the small argument asymptotics coincide with the large order asymptotics \cite[Sec. 2.1]{barnett2008stability}. See also \cite{antoine2013dilute} for related spectral and condition number estimates of the single layer operator, this time in dilute media.
\end{remark}

\begin{figure}[htbp]
	\captionsetup[subfloat]{labelformat=empty}
	\centering
	\begin{tabular}{ccc}
	\subfloat[]{\includegraphics[scale=0.6]{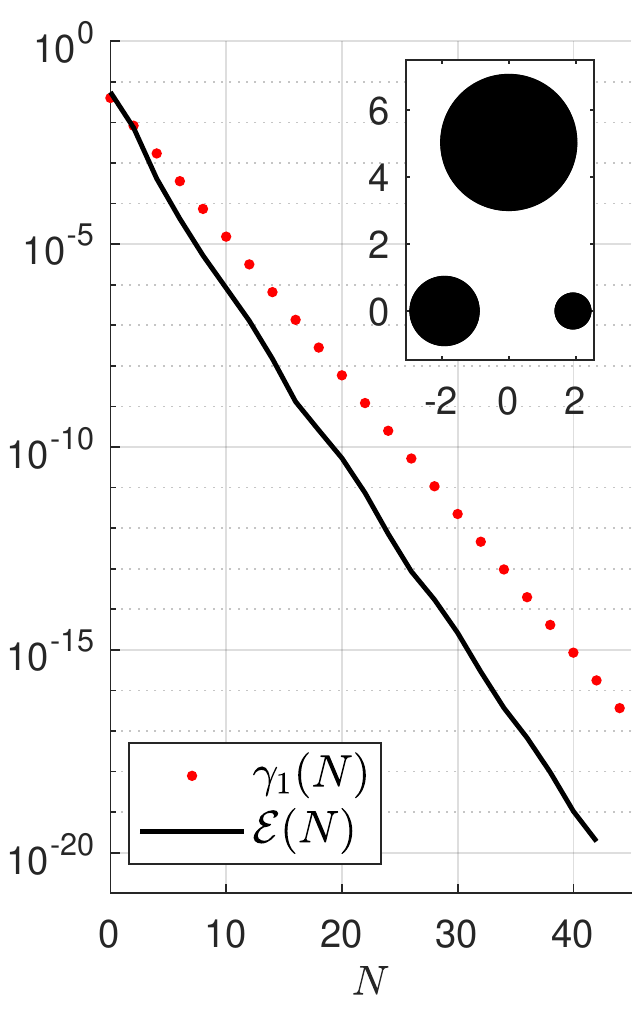}} \hspace{1.5em}
	\subfloat[]{\includegraphics[scale=0.6]{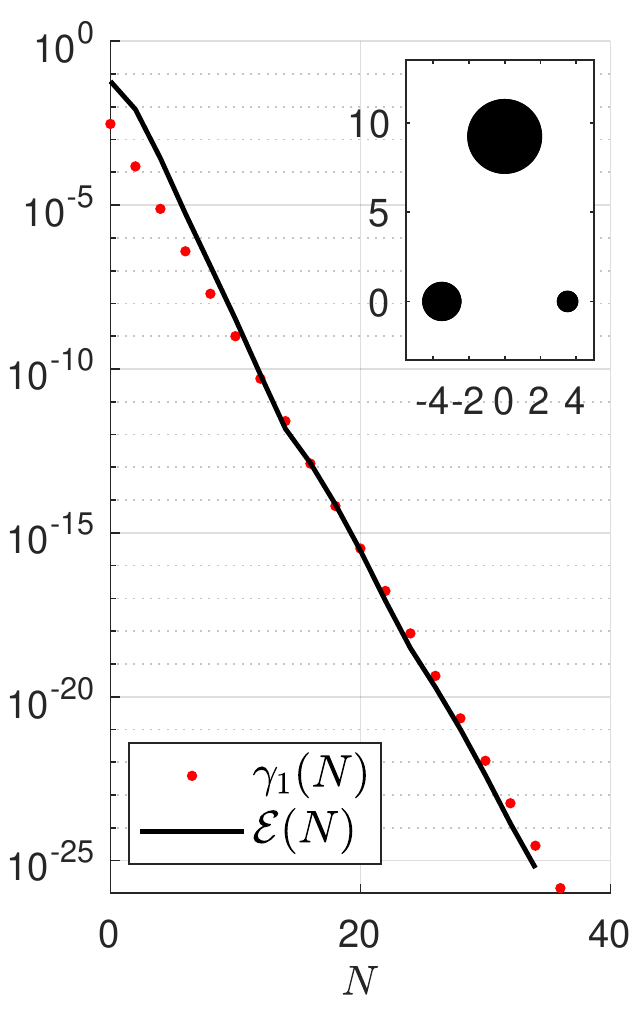}} \hspace{1.5em}
	\subfloat[]{\includegraphics[scale=0.6]{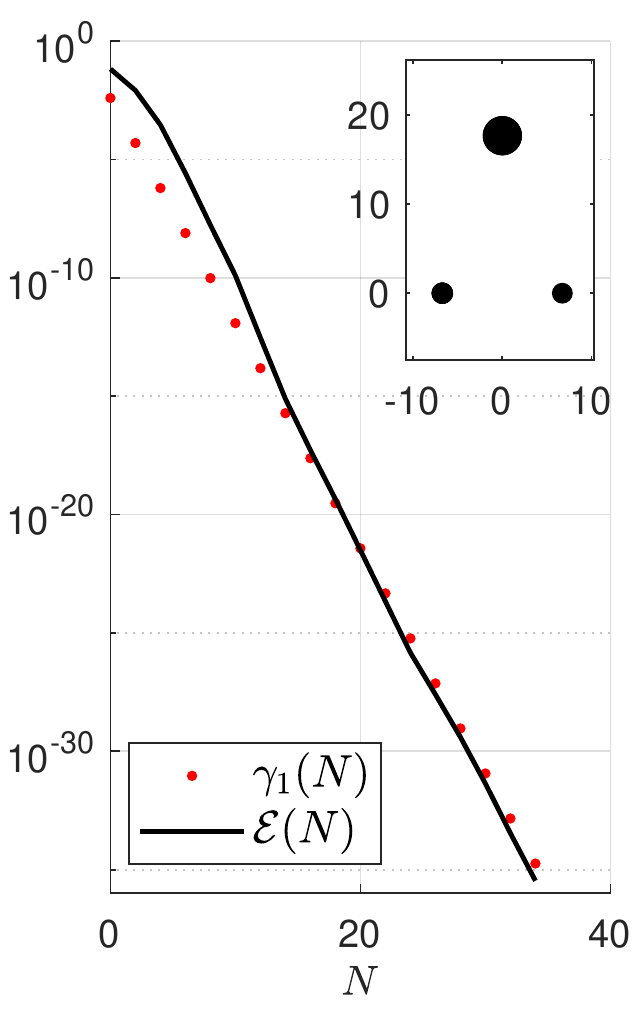}} \\
	\subfloat[]{\includegraphics[scale=0.6]{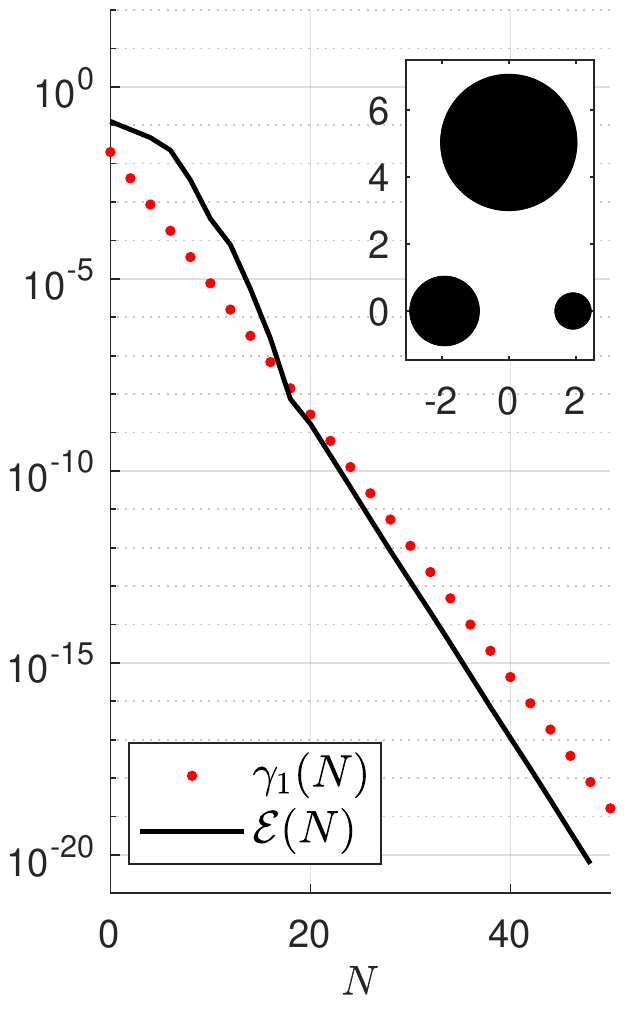}} \hspace{1.5em}
	\subfloat[]{\includegraphics[scale=0.6]{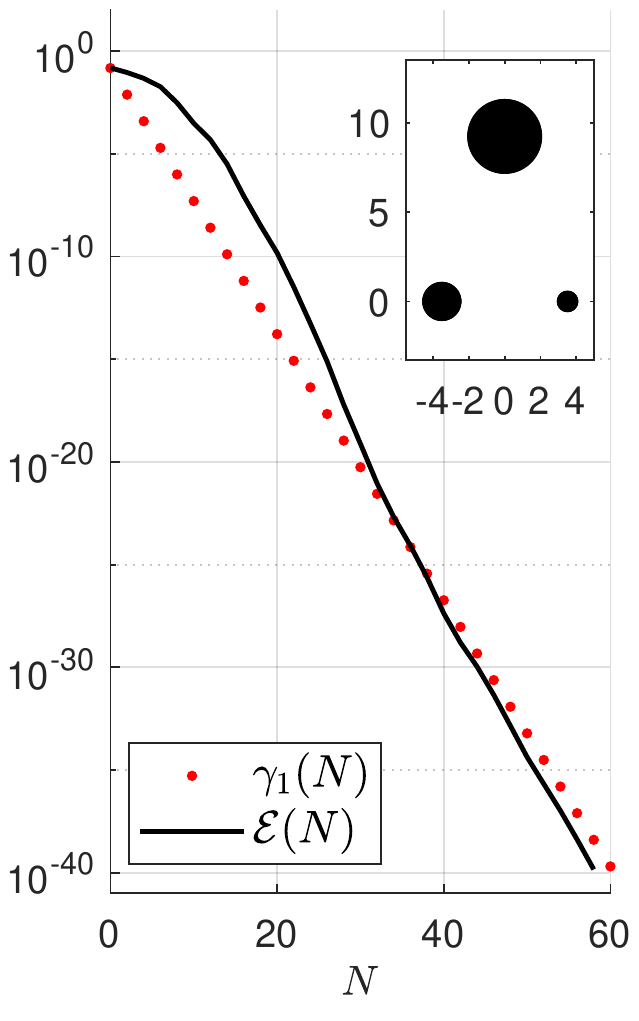}} \hspace{1.5em}
	\subfloat[]{\includegraphics[scale=0.6]{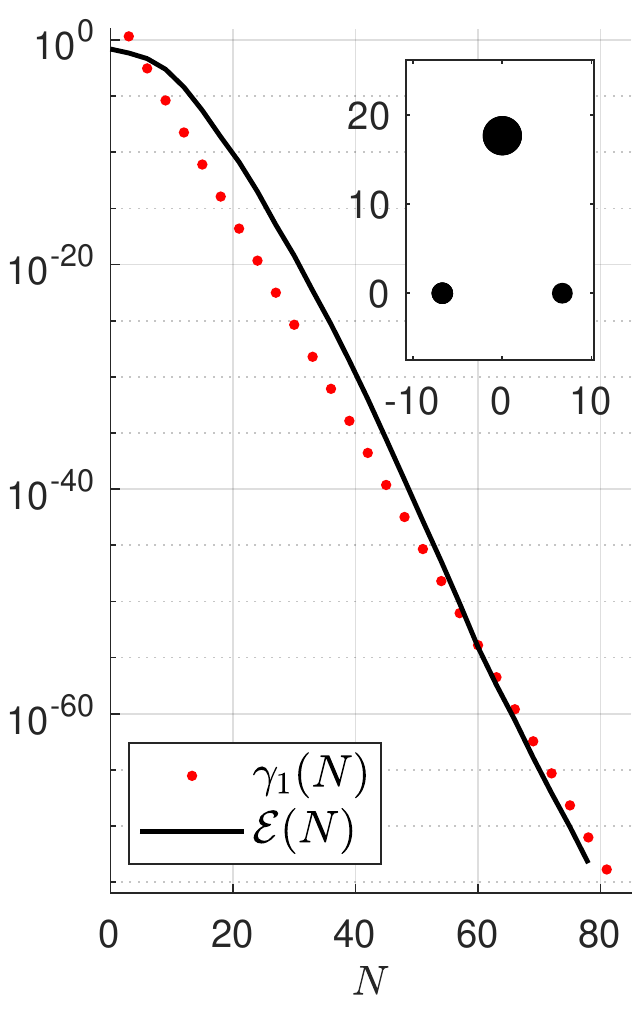}} \\
	\subfloat[]{\includegraphics[scale=0.6]{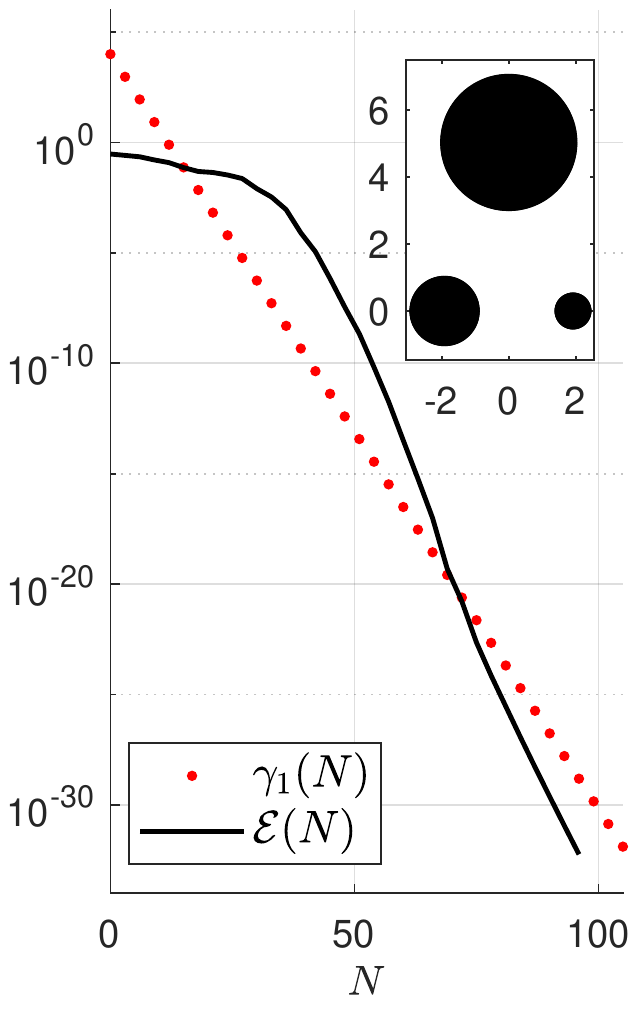}} \hspace*{\fill}
	\end{tabular}
	\caption{The convergence of the approximation error $\mathcal{E}(N)$ of the MEM, and the bound $\gamma_1(N)$ from \cref{thm:mem-app-err-bnd}, with respect to a range of geometrical configurations and wavenumbers.
	\label{fig:mem-bnds}}
\end{figure}

The results established in the previous two lemmas are succinctly summarized in the following corollary which will be used in \cref{lem:A-minus-A-N}.
\begin{corollary} \label{cor:A-pq-mn-bnd}
As $N\to \infty$, it holds that
\begin{align*}
\sumainb{m}{\bbZ_N^c} \sumz{n} |A_{mn}^{pq}|^2 \lesssim \bigg(\frac{a_p}{d_{pq}-a_q}\bigg)^{2N}, \quad \quad \sumainb{m}{\bbZ_N^c} \sumz{n} |A_{nm}^{pq}|^2 \lesssim \bigg(\frac{a_q}{d_{pq}-a_p}\bigg)^{2N}.
\end{align*}
\begin{proof}
By \cref{lem:A-pq-mn-series-bnd} and \cref{lem:sigma-series-bnd}, and noting that $0 < a_q < d_{pq}$ for $p,q\in \bbN_M$, we have
\begin{align*}
\sumainb{m}{\bbZ_N^c} \sumz{n} |A_{mn}^{pq}|^2 & \lesssim \sum_{m>N} \bigg(\frac{a_p}{d_{pq}}\bigg)^{2m} + \bigg(\frac{a_p}{d_{pq}-a_q}\bigg)^{2m} \lesssim \sum_{m>N} \bigg(\frac{a_p}{d_{pq}-a_q}\bigg)^{2m} \lesssim \bigg(\frac{a_p}{d_{pq}-a_q}\bigg)^{2N}.
\end{align*}
The result for the other series can be obtained in a similar fashion.
\end{proof}
\end{corollary}
\begin{lemma} \label{lem:A-minus-A-N}
As $N\to \infty$, it holds that
\begin{align*}
||\bbA-\tilde{\bbA}(N)||_{-1/2,-1/2} \lesssim \max_{\substack{p,q\in\bbN_M \\ q \neq p}} \bigg(\frac{a_{p}}{d_{pq}-a_{p}}\bigg)^N.
\end{align*}
\begin{proof}
Since the operator norm is dominated by the Hilbert-Schmidt norm, we have
\begin{align} \label{eq:A-pq-mn-exp}
||\bbA-\tilde{\bbA}||_{-1/2,-1/2}^2 & \le ||\bbA-\tilde{\bbA}||_{\text{HS}}^2 =  \sum_{p\in\bbN_M}\sum_{\substack{q\in\bbN_M \\ q \neq p}}\sumz{m}\sumz{n}|A_{mn}^{pq}-\tilde{A}_{mn}^{pq}|^2.
\end{align}
We decompose the two inner series as follows:
\begin{align} \label{eq:A-pq-mn-split}
\begin{split}
\sumz{m}\sumz{n}|A_{mn}^{pq}-\tilde{A}_{mn}^{pq}|^2 & = \bigg(\sum_{m\in \bbZ_N^c}\sumz{n} + \sum_{m\in \bbZ_N}\sum_{n\in \bbZ_N^c} + \sum_{m\in \bbZ_N}\sum_{n\in \bbZ_N}\bigg)|A_{mn}^{pq}|^2 \\
& = \bigg(\sum_{m\in \bbZ_N^c}\sumz{n} + \sum_{m\in \bbZ_N}\sum_{n\in \bbZ_N^c}\bigg)|A_{mn}^{pq}|^2,
\end{split}
\end{align}
since we have from \cref{eq:A-pq-mn-N} that $A_{mn}^{pq}=\tilde{A}_{mn}^{pq}$ for $m,n \in \bbZ_N$, and $\tilde{A}_{mn}^{pq} = 0$ otherwise.

Consider, for a moment, the second series on the right hand side of the above expression. It holds that
\begin{align*}
\sumainb{m}{\bbZ_N} \sum_{n\in \bbZ_N^c} |A_{mn}^{pq}|^2 = \sumainb{n}{\bbZ_N} \sum_{m\in \bbZ_N^c} |A_{nm}^{pq}|^2 =  \sum_{m\in \bbZ_N^c} \sumainb{n}{\bbZ_N}  |A_{nm}^{pq}|^2 \le \sum_{m\in \bbZ_N^c} \sumz{n} |A_{nm}^{pq}|^2.
\end{align*}
Substituting this into \cref{eq:A-pq-mn-split}, and applying \cref{cor:A-pq-mn-bnd}, we obtain
\begin{align*}
\sumz{m}\sumz{n}|A_{mn}^{pq}-\tilde{A}_{mn}^{pq}|^2
\le \sum_{m\in \bbZ_N^c}\sumz{n}\bigg(|A_{mn}^{pq}|^2 + |A_{nm}^{pq}|^2\bigg)
\lesssim \bigg(\frac{a_p}{d_{pq}-a_q}\bigg)^{2N} + \bigg(\frac{a_q}{d_{pq}-a_p}\bigg)^{2N}.
\end{align*}
Finally, substituting this result into \cref{eq:A-pq-mn-exp}, we find that
\begin{align*}
||\bbA-\tilde{\bbA}||_{-1/2,-1/2}^2 \lesssim \sum_{p\in\bbN_M}\sum_{\substack{q\in\bbN_M \\ q \neq p}} \bigg(\frac{a_p}{d_{pq}-a_q}\bigg)^{2N} + \bigg(\frac{a_q}{d_{pq}-a_p}\bigg)^{2N} \le \max_{\substack{p,q\in\bbN_M \\ q \neq p}} \bigg(\frac{a_p}{d_{pq}-a_q}\bigg)^{2N}.
\end{align*}
\end{proof}
\end{lemma}
In \cref{fig:mem-bnds}, we provide convergence plots that demonstrate the accuracy of the bound $\gamma_1(N)$ from \cref{thm:mem-app-err-bnd}, for an $M=3$ three cylinder system, with radii $(a_1,a_2,a_3) = (2,1,0.5)$, in the case of a point source incident wavefield, with the source located far away from the cylinders. The approximation error $\mathcal{E}(N)$ was computed using the MEM scattering library $\mu$-diff \cite{thierry2015mu}. For each subplot, $\mathcal{E}(N)$ is given by
\begin{align*}
\mathcal{E}(N)
= ||\Phi-\tilde \Phi(N)||_{-1/2} \approx ||\tilde \Phi(N_\text{ref})-\tilde \Phi(N)||_{-1/2},
\end{align*}
where
\begin{align*}
\tilde \Phi(N_\text{ref}) = \tilde{\mathbb{W}}^{-1}(N_\text{ref}) G(N_\text{ref}), \quad \quad \quad \quad
\tilde \Phi(N) = \tilde{\mathbb{W}}^{-1}(N) G(N),
\end{align*}
where $\tilde{\mathbb{W}}$ is defined after \cref{eq:bie-weak-form-matrix-equiv-disc}, and $N_\text{ref}$ is taken five higher than the largest truncation number $N$ used in the subplot. Due to the fast convergence of the MEM, this choice of $N_\text{ref}$ is sufficient for $\tilde \Phi(N_\text{ref})$ to accurately approximate $\Phi$. The block matrices and vectors used in $\tilde \Phi(N)$ were zero-padded to make them align correctly with the corresponding block matrices and vectors used in the reference solution $\tilde \Phi(N_\text{ref})$. As sub-exponential factors are asymptotically irrelevant, we can disregard the $(1+|m|^2)^{-1/2}$ term in the fractional Sobolev norm $l_s$ \cref{eq:l-s-norm} and perform $l_0$ based computations.

The plots in the first, second, and third columns correspond to closely spaced cylinders, moderately far apart cylinders, and far apart cylinders, respectively. With regards to wavenumbers, in each column: the first row features a wavenumber of $k=0.6$; the second row features a wavenumber of $k=3$;  the third row features a wavenumber of $k=15$.
These wavenumbers have been so chosen because they imply regimes in which the wavelength is: large with respect to the diameter of the mid-sized cylinder (first row); around the same size as the diameter of the mid-sized cylinder (second row); smaller than the diameter of the mid-sized cylinder (third row).
Hence we are analyzing the performance of the bound in a range of representative settings. Note that in the second and third columns, the last row is missing. This is because a very large truncation number $N$ is required to reach the asymptotic regime in the case of the large wavenumber $k=15$, which results in numerical precision issues.

It is clear that $\gamma_1(N)$ accurately characterizes the convergence of the approximation error $\mathcal{E}(N)$ in the moderately far apart and far apart regimes. However, it leaves something to be desired in the closely spaced regime; in this case, it is overly pessimistic. To obtain a more precise characterization of the approximation error in this regime we need to derive an estimate based on \cref{eq:app-err-bnd-0} rather than \cref{eq:app-err-bnd-1}, since the former expression takes account of the multiple scattering of the incident wavefield among the cylinders. But this expression is problematic since it features the full solution of the untruncated problem $\Phi$ on the right-hand side, a term for which a closed-form expression does not exist. To overcome this difficulty, in the next section we develop a first-order scattering approximation which allows for a more accurate characterization of the convergence of the MEM in all regimes, but particularly in the closely spaced regime in which $\gamma_1(N)$ becomes overly pessimistic.

\subsection{A first-order scattering approximation for closely spaced cylinders}
Recall that the MEM formulation of the scattering problem involves finding $\Phi$ such that
\begin{align} \label{eq:mem-prob}
(\bbI + \bbA)\Phi = G.
\end{align}
Note that if we disregard the operator $\bbA$, we simply have $\Phi = G$. This is essentially a MEM problem in which the incident field scatters off the cylinders but the subsequent multiple scattering interactions among the cylinders are neglected. As we are already using $\Phi$ to represent the full solution, let us denote by $\hat \Phi$ the solution to this single-scattering problem:
\begin{align} \label{eq:mem-prob-fir-ord}
\hat \Phi = G.
\end{align}
Denote by $\Phi_{\text{diff}} := \Phi - \hat \Phi$. This means the original solution can be decomposed as $\Phi = G + \Phi_{\text{diff}}$. Substituting this expression into \cref{eq:app-err-bnd-0} and applying the triangle inequality gives
\begin{align*} 
\cE(N) \le ||(\bbI + \tilde{\bbA}(N))^{-1}||_{-1/2,-1/2} \bigg(||G -\tilde G(N)||_{-1/2} + ||(\bbA-\tilde{\bbA}(N))G||_{-1/2} + ||(\bbA-\tilde{\bbA}(N))\Phi_{\text{diff}}||_{-1/2} \bigg).
\end{align*}
We can write this as
\begin{align} \label{eq:app-err-decomp}
\cE(N) & \le \mathcal{E}^{(1)}(N) + \mathcal{E}^{(\text{diff})}(N),
\end{align}
where
\begin{align*}
\mathcal{E}^{(1)}(N) & := ||(\bbI + \tilde{\bbA}(N))^{-1}||_{-1/2,-1/2} \ \bigg(||G -\tilde G(N)||_{-1/2} + ||(\bbA-\tilde{\bbA}(N))G||_{-1/2}\bigg), \\
\mathcal{E}^{(\text{diff})}(N) & := ||(\bbI + \tilde{\bbA}(N))^{-1}||_{-1/2,-1/2} \  ||(\bbA-\tilde{\bbA}(N))\Phi_{\text{diff}}||_{-1/2}.
\end{align*}
Here, $\mathcal{E}^{(1)}(N)$ is the component of the bound \cref{eq:app-err-decomp} associated with the first-order scattering event, while $\mathcal{E}^{(\text{diff})}(N)$ is associated with subsequent multiple scattering interactions. In some sense, the term $||(\bbA-\tilde{\bbA}(N))G||_{-1/2}$ in $\mathcal{E}^{(1)}(N)$ can be viewed as a measure of the approximation error associated with the incident wavefield striking cylinder $\Omega_p$ and being transferred once to cylinders $\Omega_q$, for $p,q \in \bbN_M$ with $q \neq p$. This is in contrast to the term $||(\bbA-\tilde{\bbA}(N))\Phi||_{-1/2}$ in \cref{eq:app-err-bnd-0} which can be viewed as a measure of the approximation error associated with the incident wavefield striking the cylinders and undergoing an infinite number of multiple reflections among them.

As long as the cylinders are not close enough together such that the effect of multiple scattering becomes comparable to the initial single scattering effect, $\mathcal{E}^{(1)}(N)$ provides the dominant contribution to the bound \cref{eq:app-err-decomp}. Hence, the bound on $\cE(N)$ in \cref{eq:app-err-bnd-0} is well approximated by $\mathcal{E}^{(1)}(N)$ which accurately characterizes the convergence. Thus, our plan now is to explicitly derive the rate of convergence of $\mathcal{E}^{(1)}(N)$. Numerical simulations, which we present later, verify that $\mathcal{E}^{(1)}(N)$ does indeed describe the convergence of $\mathcal{E}(N)$ in the closely spaced regime in which the bound from \cref{thm:mem-app-err-bnd} becomes overly pessimistic. Of course, if some cylinders are brought very close together, this first-order scattering approximation breaks down, since in that event $\mathcal{E}^{(\text{diff})}(N)$ can become quite significant. In any case, we will shortly derive the rate of convergence of $\mathcal{E}^{(1)}(N)$ for the cases of point source and plane wave incident wavefields. 

First, however, we make a brief remark about how the approximation we have just described can be connected to the method of reflections \cite{ciaramella2017review}, which is also known as the boundary decomposition method \cite{balabane2004boundary}. This approach amounts to treating the scattering of waves between objects in an iterative fashion. In fact, there are several different types of methods of reflections, the one we consider below is known as the parallel method of reflections. Denote by $\Omega_p^+ := \mathbb{R}^2 \setminus \Omega_p$ the region exterior to cylinder $\Omega_p$, for $p \in \bbN_M$. The scattered field in a Helmholtz multiple scattering problem can be decomposed as $u^\text{s} = \sum_{p\in \bbN_M} u_p^\text{s}$ where $u_p^\text{s}$ satisfies the Sommerfeld radiation condition, and
\begin{align*}
(\Delta + k^2) u_p^\text{s} = 0, \ \text{in} \ \Omega_p^+, \quad \quad u_p^\text{s} = - u^\text{inc} - \sum_{\substack{q \in \bbN_M\\ q \neq p}} u_q^\text{s}, \ \text{on} \ \Gamma_p.
\end{align*}
Furthermore, we can decompose $u_p^\text{s}$ as $u_p^\text{s} = \sum_{n=0}^\infty u_{p,n}^s$, where $u_{p,n}^s$ satisfies the Sommerfeld radiation condition, with
\begin{align} \label{eq:met-of-ref}
(\Delta + k^2) u_{p,0}^s = 0, \ \text{in} \ \Omega_p^+, \quad \quad u_{p,0}^s = - u^\text{inc}, \ \text{on} \ \Gamma_p,
\end{align}
and for $n\neq 0$,
\begin{align} \label{eq:met-of-ref-mult}
(\Delta + k^2) u_{p,n}^s = 0, \ \text{in} \ \Omega_p^+, \quad \quad u_{p,n}^s = - \sum_{\substack{q \in \bbN_M\\ q \neq p}} u_{q,n-1}^\text{s}, \ \text{on} \ \Gamma_p.
\end{align}
Now, in terms of the method of reflections, the single-scattering problem is associated with \cref{eq:met-of-ref}. In this problem only the initial single-scattering event is considered; the subsequent multiple scattering events are represented by \cref{eq:met-of-ref-mult}. So the method of reflections problem \cref{eq:met-of-ref} corresponds directly to the MEM single-scattering problem \cref{eq:mem-prob-fir-ord}.

This connection between our first-order scattering approximation and the method of reflections suggests that it may be possible to characterize the higher order multiple scattering effects in the MEM in an iterative fashion, similar to how \cref{eq:met-of-ref-mult} iteratively provides the higher order effects for the method of reflections. However, this is not entirely straightforward, as there are certain conditions that need to be met so that the method of reflections series solution converges \cite{balabane2004boundary,haibing2013decomposition}. A means of overcoming the convergence issue could be to use the averaged parallel method of reflections (APMR) introduced in \cite{laurent2017method}, which can be viewed as a relaxtion applied to the standard parallel method of reflections that leads to a convergent series. We now derive the convergence of the first-order scattering approximation $\mathcal{E}^{(1)}(N)$.

\begin{proof}[Proof of \cref{thm:mem-app-err-bnd-ps-pw}]
The proof is practically the same as that of \cref{thm:mem-app-err-bnd}, the difference being that instead of estimating $||\bbA-\tilde \bbA(N)||_{-1/2,-1/2}$, we have to estimate $||(\bbA-\tilde \bbA(N))G||_{-1/2}$. This can be accomplished in a similar manner to how $||\bbA-\tilde \bbA(N)||_{-1/2,-1/2}$ was estimated in \cref{lem:A-minus-A-N}, so we only highlight the key differences. It is straightforward to show that
\begin{align} \label{eq:A-pq-mn-exp-ps}
||(\bbA-\tilde{\bbA}(N))G||_{-1/2}^2 & \le \sum_{p\in\bbN_M}\sum_{\substack{q\in\bbN_M \\ q \neq p}}\sumz{m}\sumz{n}(1+|m|^2)^{-1/2} |A_{mn}^{pq}-\tilde{A}_{mn}^{pq}|^2 |g_n^q|^2.
\end{align}
The inner summations can decomposed as 
\begin{align} \label{eq:A-pq-mn-g-q-n-decomp-ps}
\begin{split}
\sumz{m}\sumz{n}(1+|m|^2)^{-1/2} |A_{mn}^{pq}-\tilde{A}_{mn}^{pq}|^2 |g_n^q|^2 & \le \sum_{m\in \bbZ_N^c}\sumz{n}(1+|m|^2)^{-1/2}|A_{mn}^{pq}|^2|g_n^q|^2 \\
& + \sum_{m\in \bbZ_N^c} \sumz{n}(1+|m|^2)^{-1/2} |A_{nm}^{pq}|^2|g_n^q|^2.
\end{split}
\end{align}
We can bound $\sumainb{m}{\bbZ_N^c} \sumz{n} (1+|m|^2)^{-1/2} |A_{mn}^{pq}|^2|g_n^q|^2$ as outlined in \cref{lem:A-pq-mn-series-bnd}, in the process also absorbing the algebraic factor $(1+|m|^2)^{-1/2}$, to get
\begin{align} \label{eq:A-pq-mn-g-q-n-bnd-ps}
\sumainb{m}{\bbZ_N^c} \sumz{n} (1+|m|^2)^{-1/2} |A_{mn}^{pq}|^2 |g_n^q|^2 & \lesssim \sum_{m>N} \bigg(\frac{a_p}{d_{pq}}\bigg)^{2m} + \sum_{m>N} \sum_{n=1}^\infty \sigma^{pq}(m,n), 
\end{align}
where, this time, $\sigma^{pq}(m,n)$ depends on the incident wavefield. First, we deal with the point source case, for which we have
\begin{align*}
\sigma^{pq}(m,n) = \bigg(\frac{m+n}{m}\bigg)^{2m}\bigg(\frac{m+n}{n}\bigg)^{2n}\bigg(\frac{a_p}{d_{pq}}\bigg)^{2m}\bigg(\frac{a_q^2}{d_{pq}d_{qx_0}}\bigg)^{2n}.
\end{align*}
By the same approach used in \cref{lem:sigma-series-bnd}, we find that as $m\to \infty$,
\begin{align} \label{eq:sig-g-bnd-ps}
\sum_{n=1}^\infty \sigma^{pq}(m,n) \lesssim \bigg(\frac{a_pd_{qx_0}}{d_{pq}d_{qx_0}-a_q^2}\bigg)^{2m}.
\end{align}
So, in light of \cref{eq:A-pq-mn-exp-ps}, \cref{eq:A-pq-mn-g-q-n-decomp-ps}, and \cref{eq:A-pq-mn-g-q-n-bnd-ps}, we can proceed along the sames lines as \cref{cor:A-pq-mn-bnd} to obtain that, as $N\to \infty$,
\begin{align*}
||(\bbA-\tilde{\bbA}(N))G||_{-1/2}^2 & \lesssim \max_{\substack{p,q\in\bbN_M \\ q \neq p}} \bigg(\frac{a_pd_{qx_0}}{d_{pq}d_{qx_0}-a_q^2}\bigg)^{2N}.
\end{align*}
Now, we handle the case of a plane wave incident wavefield for which we have that
\begin{align*}
\sigma^{pq}(m,n) = \bigg(\frac{m+n}{m}\bigg)^{2m}\bigg(\frac{m+n}{n}\bigg)^{2n}\bigg(\frac{a_p}{d_{pq}}\bigg)^{2m}\bigg(\frac{eka_q^2}{2d_{pq}}\bigg)^{2n}\bigg(\frac{1}{n}\bigg)^{2n}.
\end{align*}
Denote by $z^{pq} := (eka_q^2/(2d_{pq}))^2$.
Applying Stirling's approximation \cref{eq:stir-approx} to some of the terms in $\sigma^{pq}$, just as we did in \cref{lem:sigma-series-bnd}, we find that
\begin{align*}
\sum_{n=1}^\infty \sigma^{pq}(m,n)
& \lesssim \bigg(\frac{a_p}{d_{pq}}\bigg)^{2m} \sum_{n=1}^\infty \bigg(\frac{(m+n)!}{m!n!}\bigg)^2 (z^{pq})^n\bigg(\frac{1}{n}\bigg)^{2n}.
\end{align*}
This is very similar to the analogous expression in \cref{lem:sigma-series-bnd}, the difference being that we have a $(1/n)^{2n}$ term that was not present in that case. Applying Stirling's approximation to this term also, we get
\begin{align*}
\sum_{n=1}^\infty \sigma^{pq}(m,n) & \lesssim \bigg(\frac{a_p}{d_{pq}}\bigg)^{2m} \sum_{n=1}^\infty \bigg(\frac{(m+n)!}{m!n!}\bigg)^2 \frac{1}{(n!)^2 e^{2n}} (z^{pq})^{n}.
\end{align*}
Now, it transpires that
\begin{align*}
{}_{2}F_3(m+1,m+1;1,1,1;z^{pq}) = 1 + \sum_{n=1}^\infty \bigg(\frac{(m+n)!}{m!n!}\bigg)^2 \frac{1}{(n!)^2 e^{2n}} (z^{pq})^{n},
\end{align*}
so we have once again reduced the problem of bounding $\sum_{n=1}^\infty \sigma^{pq}(m,n)$ to that of finding a large argument asymptotic expansion of a hypergeometric function. In \cref{lem:hyper-geom-asy-pw}, we show that, as $x\to \infty$, it holds that
\begin{align*}
{}_{2}F_3(m+1,m+1;1,1,1;z) \lesssim \text{Exp}(4(m^2 z)^{1/4}).
\end{align*}
Therefore, as $m\to \infty$, we have
\begin{align*}
\sum_{n=1}^\infty \sigma^{pq}(m,n) & \lesssim 
\bigg(\frac{a_p}{d_{pq}}\bigg)^{2m} \text{Exp}(4(m^2 z^{pq})^{1/4}) \lesssim \bigg(\frac{a_p}{d_{pq}}\bigg)^{2m},
\end{align*}
where we absorbed the root-exponential term using \cref{eq:exp-bnd-notation}.
Once again, using \cref{eq:A-pq-mn-exp-ps}, \cref{eq:A-pq-mn-g-q-n-decomp-ps}, and \cref{eq:A-pq-mn-g-q-n-bnd-ps}, we proceed along the sames lines as \cref{cor:A-pq-mn-bnd} to obtain that, as $N\to \infty$,
\begin{align*}
||(\bbA-\tilde{\bbA}(N))G||_{-1/2} & \lesssim \sum_{p\in\bbN_M}\sum_{\substack{q\in\bbN_M \\ q \neq p}} \bigg(\frac{a_p}{d_{pq}}\bigg)^{N}.
\end{align*}
\end{proof}

\begin{figure}[htbp]
	\captionsetup[subfloat]{labelformat=empty}
	\centering
	\begin{tabular}{ccc}
	\subfloat[]{\includegraphics[scale=0.6]{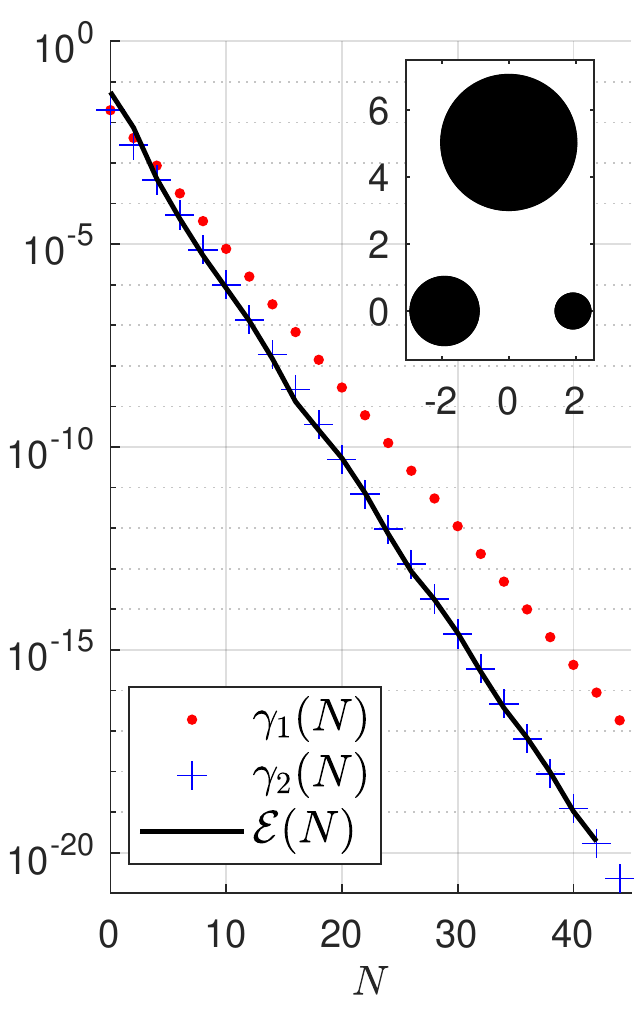}} \hspace{1.5em}
	\subfloat[]{\includegraphics[scale=0.6]{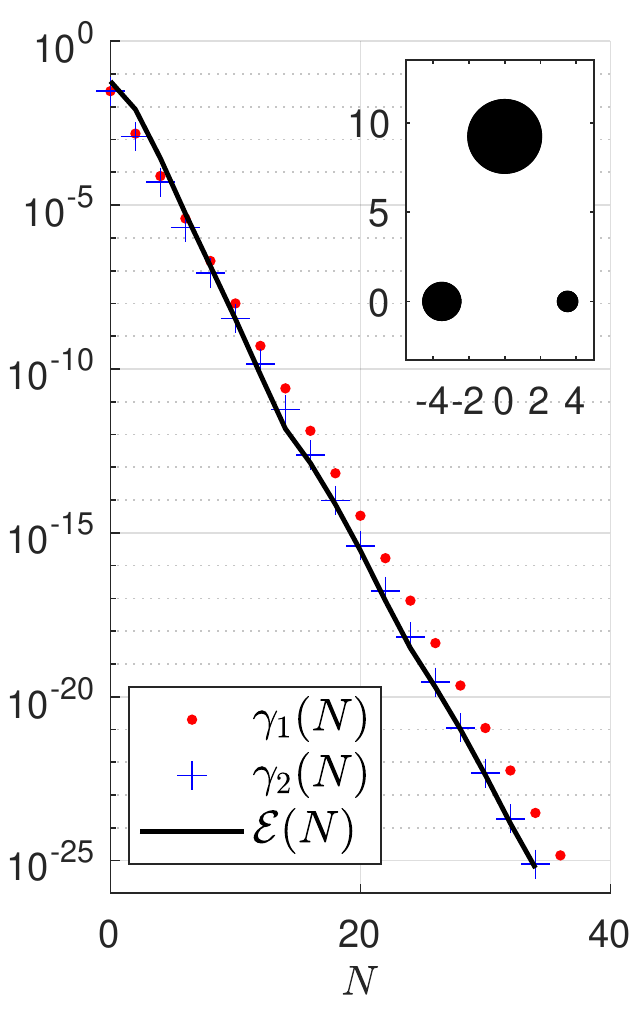}} \hspace{1.5em}
	\subfloat[]{\includegraphics[scale=0.6]{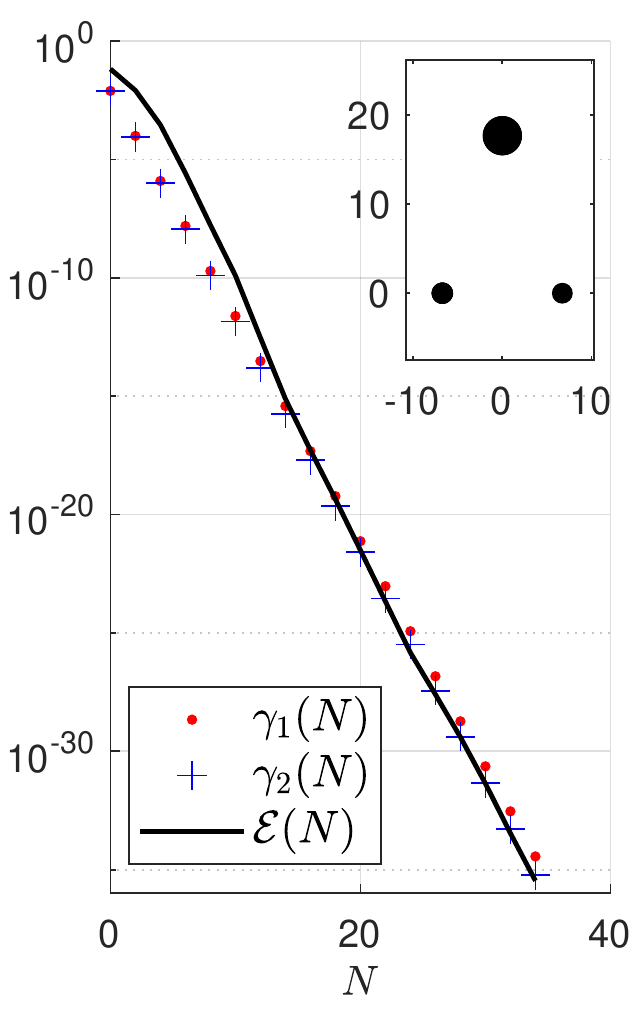}} \\
	\subfloat[]{\includegraphics[scale=0.6]{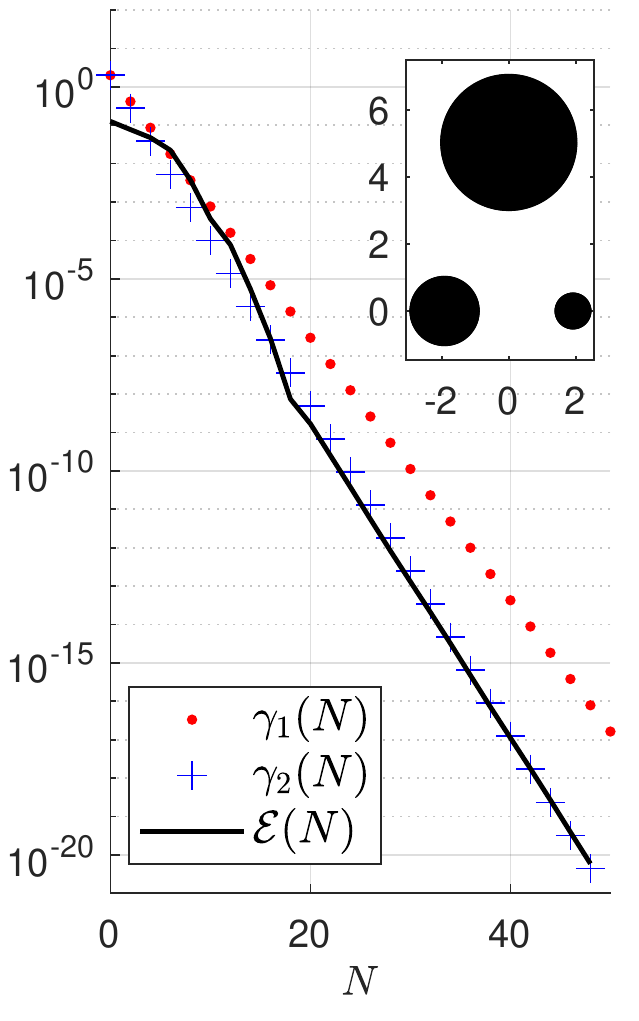}} \hspace{1.5em}
	\subfloat[]{\includegraphics[scale=0.6]{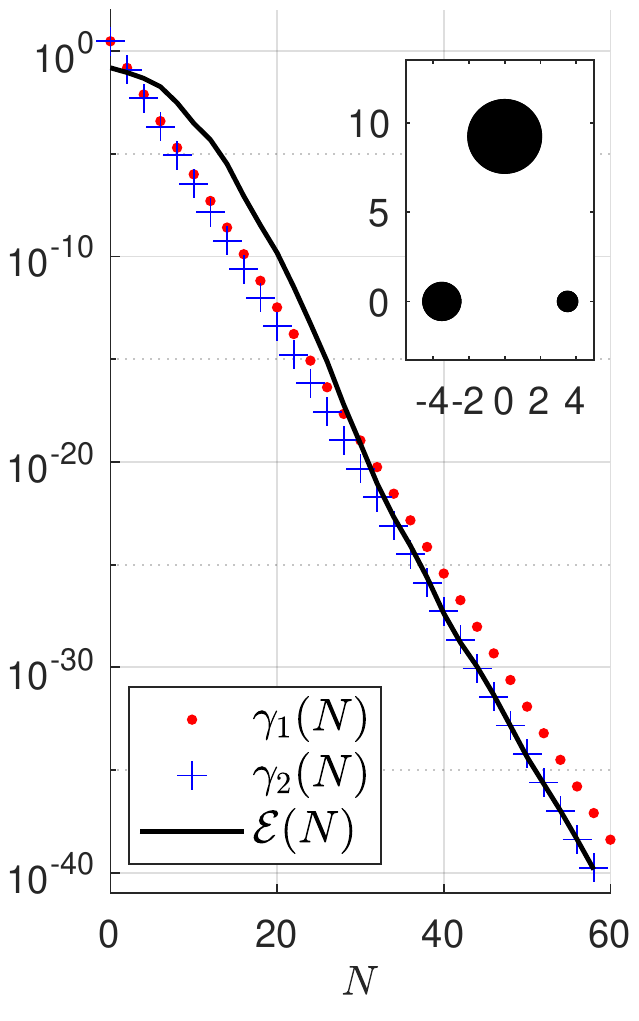}} \hspace{1.5em}
	\subfloat[]{\includegraphics[scale=0.6]{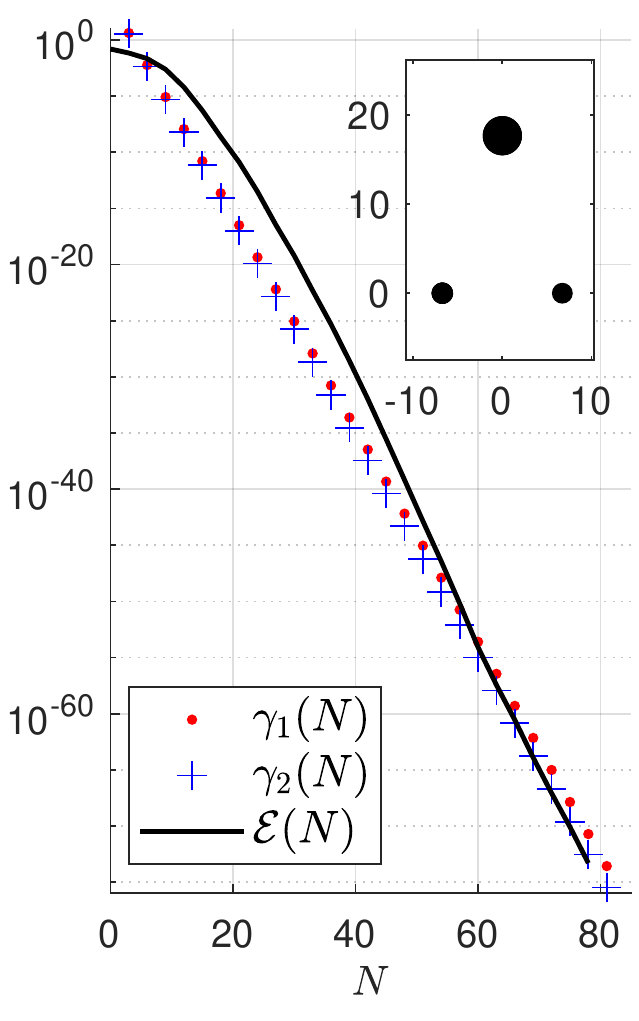}} \\
	\subfloat[]{\includegraphics[scale=0.6]{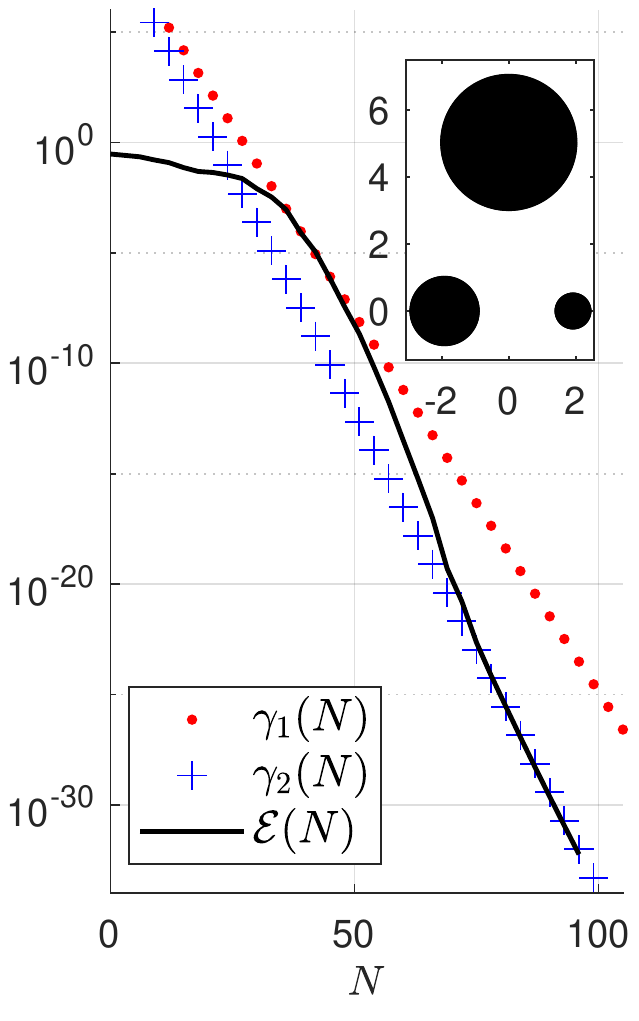}} \hspace*{\fill}
	\end{tabular}
	\caption{The convergence of the approximation error $\mathcal{E}(N)$ of the MEM, the bound $\gamma_1(N)$ from \cref{thm:mem-app-err-bnd}, and the first-order scattering bound $\gamma_2(N)$ from \cref{thm:mem-app-err-bnd-ps-pw}, with respect to a range of geometrical configurations and wavenumbers.
	\label{fig:mem-bnds-ps}}
\end{figure}

In \cref{fig:mem-bnds-ps}, we plot the convergence of the approximation error $\mathcal{E}(N)$ of the MEM, along with the bound $\gamma_2(N)$ on the first-order scattering approximation derived in \cref{thm:mem-app-err-bnd-ps-pw}, for the case of a point source located far from the cylinders; the plot for the plane wave case is very similar so we omit it. We also show the overly pessimistic bound $\gamma_1(N)$ that was derived in \cref{thm:mem-app-err-bnd}. The setting in these plots is the same as in \cref{fig:mem-bnds}, which we recall shows the convergence for the case of low (first row), medium (second row), and high (third row) wavenumbers when the cylinders are close together, a moderate distance apart, and far apart. It is clear that $\gamma_2(N)$ characterizes the convergence of the MEM better than $\gamma_1(N)$ in all regimes, with this improvement being particularly noticeable in the closely spaced regime.

Note also that $\gamma_2(N)$ may in fact slightly overestimate the convergence in the closely spaced regime, which is to be expected since in this case, while the first-order initial scattering event has a dominating effect on the convergence, the higher order multiple scattering effects, which correspond to the neglected term $\mathcal{E}^{(\text{diff})}$ in \cref{eq:app-err-decomp}, are also starting to become noticeable.

While the asymptotic rate of convergence of the MEM is wavenumber independent, it can be seen from \Cref{fig:mem-bnds-ps} that the range of values of $N$ for which $\gamma_2(N)$ accurately characterizes the convergence does have a dependence on $k$; as $k$ increases, a larger $N$ is required before the asymptotic regime is reached. This is because the functions $|J_m(ka_p)|$ and $|H_m(ka_p)|$ only begin to reach their asymptotic rates of convergence when $m \approx k a_p$. A similar phenomenon has been observed with the method of fundamental solutions \cite{barnett2008stability}, which is another meshless method featuring solution representations comprised of Bessel and Hankel functions.

An interesting question is that of when the first-order scattering approximation breaks down. Numerical investigations show that the approximation is accurate as long as the distance between the closest points on the largest cylinder and the mid-sized cylinder is greater than approximately $5/4$ times the radius of the mid-sized cylinder, irrespective of the wavenumber. Equivalently, for the wavenumbers $k = 0.6,3$, and $15$ used in \cref{fig:mem-bnds-ps}, the first-order scattering approximation is accurate as long as the distance between the closest points on the largest cylinder and the mid-sized cylinder is greater than approximately $1/8, 5/8$, and $3$ wavelengths, respectively. If the distance between the cylinders is any less than this, an approximation that takes into account higher order scattering effects would be necessary to accurately characterize the convergence of the MEM.

\section{Conclusion} \label{sec:conc}

In this work we have provided a resolution to the long-standing problem of characterizing the asymptotic rate of convergence of the approximation error of the MEM by performing a detailed convergence analysis. The system of equations we considered first arose in Row's 1955 paper \cite[Eq. (3) and Eq. (5)]{row1955theoretical}. We began by deriving a bound that is tight as long as the cylinders are not too close together. To handle the case when some cylinders are, in fact, in close proximity to one another, we formulated a first-order scattering approximation for the MEM approximation error. This approximation accounts for the initial scattering event, namely, an incident wavefield impinging on each of the cylinders which then gets reflected onto each of the other cylinders. Meanwhile, higher order repeated multiple reflections of waves among the cylinders are neglected. We derived explicit bounds on the rate of convergence of this approximation, for the cases of both point-source and plane wave incident wavefields.

While our estimates were derived based on an indirect boundary integral equation solution representation applied to a Dirichlet scattering problem, this was merely for convenience. The convergence of the MEM for other boundary conditions or boundary integral equation solution representations differs from the case we considered only by sub-exponentially increasing factors which are asymptotically irrelevant. Thus, ours is a general theory of the MEM convergence that holds for all boundary conditions and boundary integral equation solution representations.

While the primary aim of this paper was to address the long-standing question on the asymptotic convergence of the MEM, there are several avenues worthy of investigation in terms of future research. 
Firstly, one could explore the generalization of the approach outlined in this paper to the case of spheres in three dimensions. The three dimensional case is more complicated as the solution representation features not only Bessel/Hankel functions and complex exponentials, but also associated Legendre functions. Moreover, more complicated addition theorems are required \cite{martin2006multiple}. Due to extra difficulties such as these, the hypergeometric functions that could potentially arise during the analysis of the three dimensional case may turn out to be too complicated to analyze asymptotically. If this is the case, it would be interesting to investigate whether the approach employed in \cite{antoine2013dense} could be of use, as in this paper expressions somewhat similar to ours were derived without making use of hypergeometric functions.

 One could also attempt to provide a more accurate characterization of the convergence of the MEM in the case of cylinders that are almost touching; we conjectured that this may be possible by connecting our approach with a technique known as the method of reflections. Finally, since the MEM often features as a building block in other numerical methods such as, for instance, the MEM-based lattice summation techniques that arise in the field of photonic and phononic metamaterials \cite{ammari2018mathematical}, the framework we outlined in this paper could also prove helpful in obtaining rates of convergence in those approaches.

\section{Acknowledgements}
The authors wish to thank the reviewers for their insightful comments and suggestions that helped improve and clarify this manuscript.

\appendix
\begin{appendices}

\section{Asymptotic bounds for hypergeometric functions}

\begin{lemma} \label{lem:hyper-geom-asy-ps}


Let $(a/b)^2 \in(0,1)$. As $m\to \infty$, the hypergeometric function ${}_{2}F_1(m+1,m+1;1;(a/b)^2)$ is bounded as
\begin{align*}
{}_{2}F_1(m+1,m+1;1;(a/b)^2) \lesssim \bigg(\frac{b^2}{b^2-a^2}\bigg)^m \bigg(\frac{b+a}{b-a}\bigg)^m.
\end{align*}
\begin{proof}
We need the following hypergeometric function identities \cite[15.8.1,15.12.5]{olver2010nist}:
\begin{align} \label{eq:hyper-sign-change}
\begin{split}
{}_{2}F_1(a,b;c;z) & = (1-z)^{-a} {}_{2}F_1\bigg(a,c-b;c;\frac{z}{z-1}\bigg),  \\
\end{split}
\end{align}
and
\begin{align} \label{eq:hyper-asy}
\begin{split} 
{}_{2}F_1\bigg(a+\lambda,b-\lambda;c;\tfrac{1}{2}-\tfrac{1}{2}y\bigg) & = 2^{(a+b-1)/2}\frac{(y+1)^{(c-a-b-1)/2}}{(y-1)^{c/2}} \sqrt{\zeta \sinh \zeta}(\lambda+\tfrac{1}{2}a-\tfrac{1}{2}b)^{1-c} \\
& \times
\Bigg(
I_{c-1}((\lambda+\frac{1}{2}a-\frac{1}{2}b)\zeta)(1 + O(\lambda^{-2})) \\
+
& \frac{I_{c-2}(\lambda+\frac{1}{2}a-\frac{1}{2}b)\zeta}{2\lambda+a-b}
\times
\bigg((c-\tfrac{1}{2})(c-\tfrac{3}{2})(\tfrac{1}{\zeta}-\coth \zeta) \\
& +\tfrac{1}{2}(2c-a-b-1)(a+b-1)\tanh(\tfrac{1}{2}\zeta) + O(\lambda^{-2})
\bigg)
\Bigg),
\end{split}
\end{align}
where $\zeta = \text{arccosh}(y)$ and $I_c$ is the modified Bessel function of the first kind of order $c$. Substituting $a,b=m+1, c = 1$ into \cref{eq:hyper-sign-change}, we get
\begin{align} \label{eq:hyper-sign-change-explicit}
{}_{2}F_1(m+1,m+1;1;z) & = (1-z)^{-(m+1)} {}_{2}F_1\bigg(m+1,-m;1;\frac{z}{z-1}\bigg).
\end{align}
We can now apply \cref{eq:hyper-asy} to ${}_{2}F_1(m+1,-m;1;z/(z-1))$. Specifically, setting $\lambda=m,a=1,b=0$, and $c=1$ in \cref{eq:hyper-asy}, we obtain the following leading-order behavior as $m \to \infty$:
\begin{align*}
{}_{2}F_1\bigg(1+m,-m;1;\tfrac{1}{2}-\tfrac{1}{2}y\bigg) & \sim \frac{1}{\sqrt{(y+1)(y-1)}} \sqrt{\zeta \sinh(\zeta)} I_0((m+\tfrac{1}{2})\zeta).
\end{align*}
Next, setting $(1-y)/2=z/(z-1)$, which can be re-arranged as $y = (1+z)/(1-z)$,
we get
\begin{align*}
{}_{2}F_1\bigg(1+m,-m;1;\frac{z}{z-1}\bigg) & \sim \frac{1}{2\sqrt{\frac{z}{(z-1)^2}}} \sqrt{\text{arccosh}\bigg(\frac{1+z}{1-z}\bigg) \sinh\bigg(\text{arccosh}\bigg(\frac{1+z}{1-z}\bigg)\bigg)} \\
& \times I_0\bigg(\bigg(m+\frac{1}{2}\bigg)\text{arccosh}\bigg(\frac{1+z}{1-z}\bigg)\bigg).
\end{align*}
So, absorbing the algebraic factor using \cref{eq:exp-bnd-notation}, as $m \to \infty$, it holds that
\begin{align*}
{}_{2}F_1\bigg(1+m,-m;1;\frac{z}{z-1}\bigg) & \lesssim I_0\bigg(\bigg(m+\frac{1}{2}\bigg)\text{arccosh}\bigg(\frac{1+z}{1-z}\bigg)\bigg).
\end{align*}
The large argument asymptotics of the modified Bessel function \cite[9.7.1]{abramowitz1948handbook} give that $I_0(x) \lesssim e^x$ as $x\to \infty$,
and therefore, we have
\begin{align*}
{}_{2}F_1\bigg(1+m,-m;1;\frac{z}{z-1}\bigg) & \lesssim \text{Exp}\bigg({m \ \text{arccosh}\bigg(\dfrac{1+z}{1-z}\bigg)}\bigg), \quad \quad m \to \infty.
\end{align*}
The expression on the right hand side simplifies and we get
\begin{align*}
{}_{2}F_1\bigg(1+m,-m;1;\frac{z}{z-1}\bigg) & \lesssim \bigg(\frac{2}{1-\sqrt{z}} - 1\bigg)^m.
\end{align*}
Recalling \cref{eq:hyper-sign-change-explicit}, this means that, as $m \to \infty$,
\begin{align*}
{}_{2}F_1(m+1,m+1;1;z) & \lesssim \bigg(\frac{1}{1-z}\bigg)^m \bigg(\frac{2}{1-\sqrt{z}} - 1\bigg)^m.
\end{align*}
Finally, upon setting $z = (a/b)^2$, we find that as $m \to \infty$, it holds that
\begin{align*}
{}_{2}F_1(m+1,m+1;1;(a/b)^2) & \lesssim \bigg(\frac{1}{1-(a/b)^2}\bigg)^m \bigg(\frac{a+b}{b-a}\bigg)^m = \bigg(\frac{b^2}{b^2-a^2}\bigg)^m \bigg(\frac{a+b}{b-a}\bigg)^m.
\end{align*}
\end{proof}
\end{lemma}
\begin{lemma} \label{lem:hyper-geom-asy-pw}
As $m\to \infty$, the hypergeometric function ${}_{2}F_3(m+1,m+1;1,1,1;z)$ is bounded as
\begin{align*}
{}_{2}F_3(m+1,m+1;1;z) \lesssim \emph{\text{Exp}}(4(m^2 z)^{1/4}).
\end{align*}
\begin{proof}
Applying \cite[16.8.10]{olver2010nist} twice, as $m\to \infty$, it holds that ${}_{2}F_3(m,m;1,1,1;z) = {}_{0}F_3(;1,1,1;m^2z)$.
Next, applying \cite[16.11.9]{olver2010nist}, we have that as $m\to \infty$, ${}_{0}F_3(;1,1,1;m^2z) \lesssim \text{Exp}(4 (m^2 z)^{1/4})$.
\end{proof}
\end{lemma}

\end{appendices}

\bibliographystyle{siamplain}
\bibliography{bib}
\end{document}